\documentstyle[11pt]{article}
\begin{document}

\title{The $Sp_3$-grassmannian and duality
       for prime Fano threefolds of genus 9}

\author{Atanas Iliev
\thanks{Partially supported by Grant MM-1106/2001
 of the Bulgarian Foundation for Scientific Research}}
\date{}
\maketitle
\begin{abstract}
{\footnotesize
By a result of Mukai, the non-abelian Brill-Noether
locus $X = M_C(2,K:3F)$ of type II, defined by a stable
rank $2$ vector bundle $F$ of invariant $3$ over a plane
quartic curve $C$, is a prime Fano $3$-fold $X = X_{16}$
of degree $16$. The associate ruled surface $S^X = {\bf P}(F)$
is uniquely defined by $X$, and we see that for the general 
$X = X_{16}$, $S^X$ is isomorphic to the Fano surface of conics 
on $X$. The argument uses the geometry of the $Sp_3$-grassmannian 
and the double projection from a line on $X_{16}$. 
}
\end{abstract}

\bigskip
\bigskip
\bigskip

\centerline{\large \bf \S 1. Introduction}

\bigskip
\medskip

The smooth Fano 3-fold $X$ is {\it prime} if
$Pic(X) = {\bf Z}(-K_X)$, where $K_X$ is the canonical
class of $X$. The {\it degree} $d = d(X) = {-K_X}^3 = 2g(X) -2$
of a prime Fano $3$-fold $X$ is always even,
and the integer $g = g(X)$ is called the genus of $X = X_{2g-2}$.  
Prime Fano $3$-folds $X_{2g-2}$ exist \ {\it iff } 
$2 \le g = g(X) \le 12, g \not= 11$,
see \cite{I1},\cite{IP}.

Let $C$ be a smooth plane quartic curve,
and let $S = {\bf P}(F) \rightarrow C$ be a
stable ruled surface of invariant 
$e(S) = min \{ s^2 :  s \ \mbox{ a section of } S \} = 3$.
By a result of Mukai the non-abelian
Brill-Noether locus of type II
$$
X_F = M_C(2,K:3F) =
\{ E \rightarrow C: rank\ E = 2,
det\ E = det\ F \otimes K, dim \ Hom(F,E) \ge 3 \},
$$
in the moduli space
$M_C(2, det\ F \otimes K_C)$ of stable rank $2$ vector bundles on $C$
of determinant $det\ F \otimes K_C$,
is a prime Fano threefold of genus $9$,
see \cite{M3}. 
The {\it associate} Fano 3-fold
$X_S := X_F$ is uniquely defined by $S$, 
since the locus $X_F = X_{F\otimes L}$ does not depend
on the twists $F \mapsto F \otimes L$ by line bundles
$L$ on $C$, see \cite{M3}. 
In turn, the general Fano $3$-fold $X$ of genus $9$
is associate to a unique ruled surface $S = S^X$,
over the {\it associate} plane  quartic curve
$C^X$ of $X$, see {\bf (5.7)}.    

By \cite{M1}, \cite{M2}, any prime Fano 3-fold $X = X_{16}$
of genus $9$ is a linear section
of the $Sp_3$-grassmannian $\Sigma \subset {\bf P}^{13}$
by a codimension $3$ subspace
${\bf P}^{10} = {\bf P}^{10}_X \subset {\bf P}^{13}$, defined
uniquely by $X$ upto the action $\rho$ of the symplectic
group $Sp_3$ on ${\bf P}^{13}$. 

The dual projective representation $\hat{\rho}$
of $\rho$ has an invariant quartic hypersurface
$\hat{F}$ in the dual space $\hat{\bf P}^{13}$,
see {\bf (2.3)(ii)}, {\bf (2.5.a)}.   
For a given general $X = \Sigma \cap {\bf P}^{10}_X$,
the plane 
${\bf P}^2_X = ({\bf P}^{10}_X)^{\perp} \subset \hat{\bf P}^{13}$
intersects $\hat{F}$ along a smooth plane quartic curve $C_X$,
defined uniquely by $X$, see {\bf (4.1)}. 
      
The two main results in this paper are:

\

{\bf A}. \ 
{\it The Fano family ${\cal F}(X)$ of conics on the
general prime Fano $3$-fold $X$ of genus $9$
is a ruled surface over the $Sp_3$-dual plane
quartic $C_X$ of $X$}: see Theorem {\bf (4.6)}. 

\

{\bf B}. \ 
{\it The ruled surface ${\cal F}(X)$ is isomorphic to
the associate ruled surface $S^X$ of $X$; in particular
$C_X \cong C^X$}: see Theorem {\bf (6.5)}.

\

The proof of {\bf A} uses the geometry of $X$
as a subvariety of the $Sp_3$-grassmannian
$\Sigma \subset {\bf P}^{13}$.
The points $c$ of the dual quartic curve
$C_X$ of $X$ are the singular hyperplane
sections $H_c \subset \Sigma$ which contain $X$.
For any $c \in C_X$, the hyperplane section
$H_c$ has a node at the correlative pivot
point $\hat{u}(c) \in \Sigma$ of $c$,
see {\bf (2.7)} and {\bf (4.1)}.
In the lagrangian plane
${\bf P}^2_{\hat{u}(c)}$ lies a smooth
conic $q(c)$, invariant under the
action of the stabilizer subgroup
$St_{\hat{c}} \subset Sp_3$ of $H_c$,
see {\bf (3.3)}. 
The vertex surface 
$S_X = \cup_{c \in  C_X} \ q(c) \subset {\bf P}^5$
of $X$ is a ruled surface over the dual curve $C_X$
with fibers the smooth conics $q(c)$.
The points $x \in S_X$ are the vertices
of the conics $q_x \subset X$ (see {\bf (4.4)}),
and the vertex map $q_x \mapsto x$ sends the Fano surface
${\cal F}(X)$ of conics on $X$ isomorphically
onto the ruled vertex surface $S_X$ of $X$,
see {\bf (4.5.b-d)} and {\bf (4.6)}. 
This proves {\bf A}.

The proof of {\bf B} is based on the
simultaneous interpretation of the lines
$l = l_{\varepsilon} = l_L$ on the Fano
threefold $X$ as minimal sections
$C_{\varepsilon}$ of the associate surface $S^X$
and as sections $C_L$ of the vertex surface $S_X$.  
By \cite{LN}, a minimal section
$C_{\varepsilon} \subset S^X$
is represented by a space curve 
$C^l_{\varepsilon}$ of genus $3$
and degree  $7$. By \cite{M3}
the curve $C^l_{\varepsilon}$ is the
same as the curve $C^l$ defining the
inverse to a double projection
from the line $l  = l_{\varepsilon}$
on $X$, see {\bf (5.1)} and {\bf (5.4)(d)}. 
As a line on the $Sp_3$-grassmannian,
$l  = l_L$ represents a pencil of
lagrangian planes through an isotropic line
$L$, in a $3$-space ${\bf P}^3_L \subset {\bf P}^5$,
see {\bf (3.1)}.
In the intersection of ${\bf P}^3_L$
with the vertex surface $S_X$ lies a unique curve
$C_L$, which is a section of
the ruled vertex surface $S_X$.
By the $Sp_3$-geometry of the double
projection from $l \subset X \subset \Sigma$, 
the curve $C_L$ is projectively equivalent to $C^l$,
see {\bf (6.1)-(6.4)}. 
The correspondence
$C_{\varepsilon}  \leftrightarrow
l_{\varepsilon} = l = l_L \leftrightarrow C_L$
identifies the set $Min(S^X)$ of
minimal sections $C_{\varepsilon}$
of $S^X$ and the set $Min'(S_X)$
of sections $C_L$ of $S_X$.
At the end, Theorem {\bf (6.5)} shows
that the isomorphism $Min(S^X) \cong Min'(S_X)$
yields an isomorphism between the  associate
surface  $S^X$ and the isomorphic vertex
model $S_X$ of ${\cal F}(X)$. This proves {\bf B}.

\bigskip
\bigskip
\bigskip

\centerline{\large \bf \S 2.
     The ${\bf Sp_3}$-grassmannian ${\bf \Sigma \subset {\bf P}^{13}}$}

\bigskip
\medskip

{\bf (2.1)}
Let $V_6 = {\bf C}^6$ be a complex vector $6$-space,
let $\hat{V}_6 = Hom(V_6,{\bf C})$ be its dual space, 
and fix a rank six $2$-form $\alpha \in \wedge^2 \ \hat{V}_6$.
The 21-dimensional symplectic group
$Sp_3 = Sp_3^{\alpha} \subset GL_6$ is the set of
all $A \in GL_6$ which preserve the non-degenerate skew-symmetric product
$\alpha : {V}_6 \times {V}_6 \rightarrow {\bf C}$ defined by $\alpha$.
In coordinates $(e_i,x_i) = (e_1,...e_6;x_1, ..., x_6)$ on $V_6$ 
$$
Sp_3 = \{A \in GL_6:{}^tAJA = J \}
$$
where $J_{ij} = {\alpha}(e_i,e_j)$
is the skew-symmetric Gramm matrix of $\alpha$.
One can always choose the coordinates $(e_i,x_i)$
on $V_6$ such that $\alpha = x_{14} + x_{25} + x_{36}$; 
and then $J_{ij} = J_{3+i,3+j} = 0$,
$J_{i,3+j} = \delta_{ij} = -J_{3+i,j}$, $1 \le i,j \le 3$. 
%

The subspace $U \subset V_6$ is called {\it isotropic} 
if $\alpha|_U \equiv 0$, 
i.e. $\alpha (u_1,u_2)  = 0$ for any $u_1,u_2 \in U$;
and then the projective subspace
${\bf P}(U) \subset  {\bf P}(V_6) = {\bf P}^5$
is also called isotropic.

Since $\alpha$ is skew-symmetric then
any point $x \in {\bf P}^5$ is isotropic.
Since $\alpha$ is non-degenerate
then the dimension of an isotropic $U \subset V$
can't be more than $3$.  The isotropic
subspaces $U \subset V_6$ of dimension $3$,
as well their projective planes ${\bf P}(U) \subset {\bf P}^5$, 
are called {\it lagrangian}.

The {\it isotropic grassmannian} $LG_2 \subset G(2,V_6)$
is the set of all isotropic lines $L \subset {\bf P}^5$.
One can see that the line $L = {\bf P}(U_2)$
is isotropic \ {\it iff} \ its Pl\"ucker image
${\bf P}(\wedge^2 \ U_2)$ lies in the {\it isotropic}
hyperplane
$(\alpha = 0) \subset {\bf P}(\wedge^2 \ V_6)$, i.e.
$$
LG_2 = G(2,6) \cap (\alpha = 0) =
G(2,6) \cap (x_{14}+x_{25}+x_{36}=0);
$$
and since $\alpha$ is non-degenerate, $LG_2$ is a
smooth hyperplane section of $G(2,6)$.

The lagrangian grassmannian,
or the {\it $Sp_3$-grassmannian} $\Sigma = LG_3 \subset G(3,V_6)$
is the set of all lagrangian subspaces $U \subset V_6$.
The $Sp_3$-grassmannian, as well the isotropic grassmannian 
$LG_2$, is a smooth homogeneous variety of the group $Sp_3$. 
In the paper we shall use the explicit coordinate description 
of $\Sigma$ which we state below. 

For the $3$-space $U \subset V_6$, the condition  $\alpha|_U \equiv 0$  
is equivalent to the requirement that its Pl\"ucker image $pl(U)$ = 
${\bf P}(\wedge^3 \ U) \in {\bf P}^{19} = {\bf P}(\wedge^3 \ V_6)$
lies in the codimension $6$ subspace
${\bf P}^{13} = {\bf P}(V_{14}) = (\hat{V}_6 \wedge \alpha = 0)
\subset {\bf P}^{19}, \ \ \mbox{ i.e.}$ 

\

\centerline{$\Sigma = G(3,6) \cap (\hat{V}_6 \wedge \alpha = 0) 
= G(3,6) \cap (x_{i14} + x_{i25} + x_{i36} = 0, 1 \le i \le 6)$.}   

\

Let $U_o = <e_1,e_2,e_3>$ and $U_{\infty} = <e_4,e_5,e_6>$. 
Then $V_6 = U_o \oplus U_{\infty}$, and in homogeneous 
coordinates 
$$
(u:X:Y:z) = (u:(x_{ij}):(y_{ij}):z) :=
\big( x_{123} :
\left( \begin{array}{ccc}
x_{423} & x_{143} & x_{124} \\
x_{523} & x_{153} & x_{125} \\
x_{623} & x_{163} & x_{126}
\end{array} \right) :
\left( \begin{array}{ccc}
x_{156} & x_{416} & x_{451} \\
x_{256} & x_{426} & x_{452} \\
x_{356} & x_{436} & x_{453}
\end{array} \right) :
x_{456} \big)
$$
on ${\bf P}^{19} = {\bf P}(\wedge^3(U_o \oplus U_{\infty})) 
= {\bf P}(\wedge^3 U_o \oplus \wedge^2 U_o \otimes U_{\infty} 
  \oplus U_o \otimes \wedge^2 U_{\infty} \oplus \wedge^3 U_{\infty})$, 
the equations $\hat{V}_6 \wedge \alpha = 0$ of 
${\bf P}^{13} \subset {\bf P}^{19}(u:X:Y:z)$ become 
${}^tX = X, \ {}^tY = Y$; i.e. 
$$
{\bf P}^{13} = {\bf P}^{13}(u:X:Y:z), \  {}^tX = X, \ {}^tY = Y.
$$  

Let $G(3,6)^o = G(3,6) - H^G_o$
be the complement to the Schubert hyperplane section
$H^G_o = {\sigma}_{100}(U_o) = \{ U \in G(3,6):
dim(U \cap U_o) > 0  \} = G(3,6) \cap (u = 0)$. 
Let also $H_o \subset \Sigma$ be the hyperplane section 
$\Sigma \cap H^G_o = \Sigma \cap (u = 0)$, and 
$\Sigma^o = \Sigma - H_o$. 
The $3$-spaces $U = U_X \in G(3,6)^o$ are parameterized
$1:1$ by the linear maps $X_U: U_o \rightarrow U_{\infty}$,
by
$$
U = <e_1 + X(e_1),e_2 + X(e_2),e_3+X(e_3)>
\ \longmapsto \ X \in Hom(U_o,U_{\infty}) \cong \otimes^2 {\bf C}^3,  
$$
and then $pl(U) = (1: X : \wedge^2 \ X : \det \ X)$. 
Moreover, the $3$-space $U = U_X \in G(3,6)^o$ is lagrangian  
(i.e. $U \in \Sigma^o$) \ {\it iff} \ the matrix $X = X_U$ is symmetric. 
Therefore $\Sigma$ is the natural projective
compactification of the isomorphic image
$\Sigma^o  = \Sigma \cap (u \not= 0)$ of the affine space 
${\bf C}^6 =  Sym^2 \ {\bf C}^3$ under the $exp$-map:
$$
exp: Sym^2 \ {\bf C}^3 \rightarrow  {\bf P}^{13}(u:X:Y:z), \
\ X \mapsto exp(X) = (1: X : \wedge^2 \ X : \det \ X).
$$

As a subvariety of ${\bf P}^{13}(u:X:Y:z)$, ${}^tX = X$, ${}^tY = Y$,
the smooth $6$-fold $\Sigma$ is defined by the 
$21$ projectivized quadratic Cramer equations

\

$({\bf \ast})$ \hspace{3.5cm}
$\wedge^2 \ X = uY, \ \wedge^2 \ Y = zX, \ XY = uzI_3$

\

\noindent (where $I_3$ is the unit $3 \times 3$ matrix), 
which follow from the local $exp$-parameterization 
of $\Sigma$ by the skew-powers of symmetric
$3 \times 3$ matrices.

\bigskip

{\bf (2.2)} {\bf Remark}. \ 
The next Theorem {\bf (2.3)} and Lemma {\bf (2.7)} below
are particular cases of the more general Propositions 5.10-5.11
and Proposition 8.2 in \cite{LM1}, where are described the orbits
in the enveloping spaces and the singular hyperplane sections
of the four varieties $\Sigma_{A} = G_{\omega}(A^3,A^6)$,
in the third row of the Freudenthal magic square.
The four varieties $G_{\omega}({A}^3,{A}^6)$
correspond to the four complex composition algebras 
${A} = {\bf A}_{\bf R} \times {\bf C}$,
where ${\bf A}_{\bf R} = {\bf R}, {\bf C}, {\bf H}$
and ${\bf O}$ = the real numbers, the complex numbers,
the quaternions and the octonions.   
The $Sp_3$-grassmannian $\Sigma = \Sigma_R$
corresponds to the case ${A} = {R}$;  
see also \S 3 in \cite{D}, where these results are proved
in case $A = {C}$, where $\Sigma_C = G(3,6)$. 

\

The subspace $V_{14} \subset \wedge^3 V_6$
is the irreducible representation
space of $Sp_3$, defined by the weight $\wedge_3$, \ 
see \cite{FH}, p. 258. 
We shall denote by $\rho$ the induced projective action
of $\wedge_3$ on ${\bf P}^{13} = {\bf P}(V_{14})$.  
%

\

{\bf (2.3)} \ {\bf Theorem of Segre for ${\bf Sp_3}$}.
\ {\sl
The action
$\rho : Sp_3 \times {\bf P}^{13} \rightarrow {\bf P}^{13}$
has $4$ orbits

\

\centerline{${\bf P}^{13} = {\Sigma} \cup ({\Omega} - {\Sigma})
\cup ({F} - {\Omega}) \cup ({\bf P}^{13} - F)$,  \ \ and:}

\

{\bf (i)} \ $\Sigma$ is the isomorphic Pl\"ucker
image of the $Sp_3$-grassmannian $LG_3$.

{\bf (ii)} \ $F \subset {\bf P}^{13}$ \ is a quartic hypersurface.
Moreover \ $F = {\cup}_{u \in {\Sigma}} \ {\bf P}^6_u$,
where ${\bf P}^6_u$ denotes the tangent
projective $6$-space at the point $u \in \Sigma$.

{\bf (iii)} \ $x \in {\bf P}^{13} - \Omega$ \ $\Leftrightarrow$ \ 
there exists a unique secant or tangent line $l_x$ to $\Sigma$,
s.t. $x \in l_x$. Moreover 
$x \in F - \Omega$ \ $\Leftrightarrow$ \ $l_x$ is tangent to $\Sigma$. 
}

\

For the point $u \in \Sigma$ we shall denote by
${\bf P}^2_u \subset {\bf P}^5$ the lagrangian
plane of $u$.  

\
 
{\bf (2.4)} \ {\bf Lemma}.
\ {\sl
Let $u \in \Sigma$,
and let ${\bf P}^6_u$ be the tangent projective
space of $\Sigma \subset {\bf P}^{13}$ at $u$.
Then ${\bf P}^6_u \subset F$, and:  
\ {\bf (i)}  
$K_u = \{ v\in\Sigma : dim({\bf P}^2_v \cap {\bf P}^2_u) \ge 1\}$
= ${\Sigma} \cap {\bf P}^6_u$
is a cone with vertex $u$ over the Veronese surface. 
\ {\bf (ii)} 
$D_u = {\Omega} \cap {\bf P}^6_u$ is a cone with vertex $u$
over the symmetric determinantal cubic. 
}

\bigskip

{\bf Proof}. \ 
%
By the $Sp_3$-homogeneity of $\Sigma$ one can let  
$u = (1:0:0:0)$ in ${\bf P}^{13}(u:X:Y:z)$.
In the notation of {\bf (2.1)}, the lagrangian plane
${\bf P}^2_u = {\bf P}(U_o)$, and the projective
tangent space ${\bf P}^6_u =
{\bf P}(\wedge^3 \ U_o \oplus \wedge^2 U_o \otimes U_{\infty})$
= ${\bf P}^6_o(u:X:0:0)$.
Therefore $v \in {\Sigma} \cap {\bf P}^6$
\ {\it iff } \ $v \in  \Sigma$ and
$dim ({\bf P}^2_v \cap {\bf P}^2_u) \ge 1$.
By {\bf (2.1)}(${\bf \ast}$), 
the intersection ${\Sigma} \cap {\bf P}^6_u(u:X:0:0)$
is defined by the system $rank \ X = 1$, ${}^tX =  X$,
which implies {\bf (i)}. By {\bf (2.3)(iii)} the
intersection ${\Omega} \cap {\bf P}^6_u(u:X:0:0)$
is exactly the determinantal cubic cone
$D_u = (det \ X = 0)$ which is the locus of points
$\omega  \in {\bf P}^6_u$ through which pass
more than one secant line to $K_u$, 
see also \S 9 of \cite{K}.
{\bf q.e.d.}

\bigskip

{\bf (2.5.a)} {\bf The varieties
${\bf \Sigma \subset \Omega \subset F \subset {\bf P}^{13}}$}. \
The smooth variety $\Sigma \subset {\bf P}^{13}$
is a Fano $6$-fold of degree $16$ and of index $4$,
i.e. $K_{\Sigma} = {\cal O}_{\Sigma}(-4)$,
see e.g. \cite{M2}, \cite{M4}. 
In coordinates $(u:X:Y:z)$ as above, the $Sp_3$-invariant
quartic hypersurface $F$ is defined by the equation
$$
F(u:X:Y:z) = (uz - tr \ XY)^2 + 4u \ det \  Y + 4z \ det \ X
- 4{\Sigma}_{ij} \ det(X_{ij}).det(Y_{ij}) =  0,
$$
see \cite{KS} p. 83, or \S 5.2 and Proposition 5.8 in \cite{LM1}. 
The $9$-fold $\Omega = Sing \ F = (\nabla \ F = 0)$
is the common zero-locus of the 14 partial cubic derivatives
of $F$. It has degree $21$, and $Sing \ \Omega = \Sigma$,
see \S 2 of \cite{IR}.

\bigskip

{\bf (2.5.b)} {\bf The pivot-map}. \ 
By analogy with \S 3 in \cite{D},
call the {\it axis} of $x \in {\bf P}^{13} - \Omega$
the unique secant or tangent line $l_x$ to $\Sigma$ through $x$,
and the {\it pivots} of $x$ 
the two intersection points $u$ and $v$ of $l_x$
with $\Sigma$. 
If $x \in F - \Omega$ call the point $u(x) = u = v \in \Sigma$
the double {\it pivot} (or simply -- the pivot) of $x$.
This way, there exists a well-defined
{\it pivot}-map:
$$
piv: F - \Omega \rightarrow \Sigma , \ x \mapsto
\mbox{ the double pivot } u = u(x) \mbox{ of } x .
$$

The fiber of the regular pivot map
$piv: F - \Omega \rightarrow \Sigma$
over the point $u \in \Sigma$ coincides with
${\bf P}^6_u - D_u$, where $D_u = \Omega \cap {\bf P}^6_u$
is the determinantal cubic cone in ${\bf P}^6_u$,
see {\bf (2.3)-(2.4)}. 

\bigskip

{\bf (2.5.c)} {\bf The correlation ${\bf J_{\alpha}}$}. \ 
The dual representation of $\wedge_3$ is equivalent to $\wedge_3$.
Therefore its  projective action $\hat{\rho}$ splits
the dual projective space $\hat{\bf P}^{13} = {\bf P}(\hat{V}_{14})$
into a union of orbits $\hat{\Sigma}$, $\hat{\Omega} - \hat{\Sigma}$,
$\hat{F} - \hat{\Omega}$ and $\hat{\bf P}^{13} - \hat{F}$ projectively
equivalent to the orbits of $\rho$.
The non-degenerate form
$\alpha \in \wedge^2 \ \hat{V}_6$ defines a linear isomorphism  
$$
J = J_{\alpha} : V_6 \rightarrow \hat{V}_6 = Hom(V_6,{\bf C}),
\
J : x \mapsto \alpha(x,.).
$$

In canonical coordinates $(e_i,x_i)$, in which
$\alpha = x_{14} + x_{25} + x_{36}$, one has: \  
$J(e_1) = x_4$, $J(e_2) = x_5$, $J(e_3) = x_6$,
$J(e_4) = -x_1$, $J(e_5) = -x_2$ and $J(e_6) = -x_3$. 

The linear isomorphism $J$ induces
a projective-linear isomorphism, or a {\it correlation} 
$$
J: {\bf P}^{13} \rightarrow \hat{\bf P}^{13}, \
J: x \mapsto \hat{x} := J(x),
$$
which identifies the elements $x$ of any of the orbits
of $\rho$ with their {\it correlative} elements
$J(x) = \hat{x}$ in its corresponding orbit of $\hat{\rho}$.
The correlation $J$ commutes with taking pivots; 
in particular, for $x \in \hat{F} - \hat{\Omega}$
the correlative point $\hat{u(x)} \in \hat{\Sigma}$
of its (double) pivot $u(x) = piv(x) \in \Sigma$
is the pivot $\hat{u}(x) := piv(\hat{x})$
of its correlative point $\hat{x} \in F - \Omega$.
%

\bigskip

{\bf (2.6)} {\bf Lemma}. \
{\sl
In coordinates $(u:X:Y:z)$, the correlative pivot map

\medskip

\centerline{$\widehat{piv}: \hat{F} - \hat{\Omega}
\rightarrow \Sigma , \ x \mapsto \widehat{u(x)},$}

\medskip

\noindent
where $u(x) = piv(x)$,
coincides with the gradient map $\nabla = \nabla_F$,
defined by the $14$ cubic derivatives of
the quartic form $F(u:X:Y:z)$.
}

\bigskip

{\bf Proof}. \ 
In coordinates $(u:X:Y:z)$, let $u = (0:0:0:1) = (x_{456})$.
By the $Sp_3$-transitivity on the orbit $\hat{F}-\hat{\Omega}$,
it is enough to prove that the gradient map $\nabla = \nabla_F$
sends any $x = (0:0:Y:z) \in {\bf P}^6_u - D_u$ to the point
$\hat{u} = \widehat{piv(u)} = (\widehat{x_{456}}) = (-e_{123})
= (e_{123}) = (1:0:0:0)$.  
This is direct: the straightforward check shows that
$$
\nabla_F: (0:0:Y:z) \mapsto (4det \ Y:0:0:0) = (1:0:0:0)
$$
since $det \ Y \not= 0$ on ${\bf P}^6_u - D_u$,
see {\bf (2.3)-(2.4)}. 

\bigskip

{\bf (2.7)} \ {\bf Lemma}. \
{\sl
{\bf (i)} 
If $x \in \hat{F} - \hat{\Omega}$
then the hyperplane section $H_x \subset \Sigma$
has a node at the correlative double pivot
$\hat{u} := \hat{u(x)}$ of $x$.
Moreover $H_x$ contains the $3$-fold cone
$K_{\hat{u}} = {\bf P}^6_{\hat{u}} \cap \Sigma$. \ 
{\bf (ii)} 
For $u \in \hat{\Sigma}$ the hyperplane section
$H_u \subset \Sigma$ is the same as the Schubert
hyperplane section
${\sigma}_{100}({\bf P}^2_{\hat{u}}) \cap {\bf P}^{13} =
\{ v\in\Sigma :{\bf P}^2_v\cap{\bf P}^2_{\hat{u}}\not= \emptyset\}$.
The singular locus
$Sing \ H_u = K_{\hat{u}} = {\bf P}^6_{\hat{u}} \cap \Sigma$
is a cone, with vertex $\hat{u}$, over the Veronese surface.
}

\bigskip

{\bf Proof}: See {\bf (2.2)}.

\bigskip 
\bigskip
\bigskip

\centerline{\large \bf \S 3.
        Lines and quadrics on the ${\bf Sp_3}$-grassmannian}

\bigskip
\medskip

{\bf (3.1)} {\bf Lemma}. \ 
{\sl
{\bf (i)} \ 
Let $L \subset {\bf P}^5$ be an isotropic line. Then the set

\medskip

\centerline{$
l_L = \{ u \in \Sigma : L \subset {\bf  P}^2_u \} \subset \Sigma
$}

\medskip

\noindent
is a line, and any line in $\Sigma$ is obtained this way. 

{\bf (ii)} \ 
If ${\bf P}^3_L \subset {\bf P}^5$ is the 3-space such that

\medskip

\centerline{$
l_L = {\sigma}_{332}(L,{\bf  P}^3_L)
= \{ u \in  \Sigma \subset G(3,6):
L \subset {\bf P}^2_u \subset {\bf P}^3_L \}
$,}

\medskip

\noindent
then ${\bf P}^3_L$ does not contain other lagrangian
planes except the planes ${\bf P}^2_u, u \in l_L$.
}

\

For the line
$l = l_L = {\sigma}_{322}(L,{\bf  P}^3_L) \subset \Sigma$,
we call the isotropic line $L$ and the space ${\bf P}^3_L$
respectively the {\it axis} and the {\it space} of $l$.

\bigskip

{\bf Proof}. \ For {\bf (i)}, see e.g. \S 2.7.2 in \cite{LM2}.
In particular, $Sp_3$ acts transitively on the
family of lines $l = l_L \subset \Sigma$, since
the isotropic grassmannian $LG_2$ is a homogeneous
variety of $Sp_3$. Therefore it is enough to prove
{\bf (ii)} for a particular line $l$ on $\Sigma$. 
In coordinates $(u:X:Y:z)$, for the isotropic line
$L = {\bf P}(<e_2,e_3>) \subset {\bf P}^5$ 
the space ${\bf P}^3_L = {\bf P}(U_4^L)$
= ${\bf P}(<e_1,e_2,e_3,e_4>)$
and the line $l = l_L = <e_{123},e_{423}>$.
Therefore the Pl\"ucker image $u$ of a plane 
${\bf P}^2_u \subset {\bf P}^3_L$ will lie
in the $3$-space ${\bf P}^3_l = {\bf P}(\wedge^3 \ U_4^L)$
= ${\bf P}^3(x_{123}:x_{423}:x_{143}:x_{124})$
= ${\bf P}^3(u:x_{11}:x_{12}:x_{13})$.
By the equations {\bf (2.1)}(${\bf \ast}$),
the plane ${\bf P}^2_u$ will be lagrangian
\ {\it iff} \
$u \in {\bf P}^3_l \cap ({}^tX = X, {}^tY = Y)$
= ${\bf P}^1(u:x_{11})$ = ${\bf P}^1(x_{123}:x_{423}) = l$. 
{\bf q.e.d.}

\bigskip

{\bf (3.2)} {\bf Lemma}. \ 
{\sl
Let $x \in {\bf P}^5$. Then the set

\medskip

\centerline{$
Q_x := {\sigma}_{300}(x) \cap {\bf P}^{13} 
= \{ u \in \Sigma : x \in  {\bf P}^2_u \} \subset \Sigma
$}

\medskip

\noindent
is a smooth 3-fold quadric, and any 3-fold quadric
on  $\Sigma$ is one of the quadrics  $Q_x,  x  \in {\bf P}^5$.
}

\

For the quadric $Q = Q_x \subset \Sigma$,
we call the point $ver(Q) := x \in {\bf P}^5$
the {\it vertex} of $Q$.

\bigskip

{\bf Proof}. \ See \S 2.7.2 in \cite{LM2}; 
see also the proof of {\bf (3.3)}{\bf (i)} below. 

%

\bigskip

{\bf (3.3)} 
{\bf Proposition
: quadrics on singular hyperplane sections of ${\bf \Sigma}$}. \ 
{\sl

{\bf (i)}
Let $Q  = Q_x  \subset \Sigma$
be a $3$-fold quadric on $\Sigma$ (see {\bf (3.2)}),   
and let $H_c, c \in \hat{\bf P}^{13}$
be a hyperplane section of $\Sigma$ containing $Q$.
Then $c \in \hat{F}$;
and if $c \not\in \hat{\Omega}$, 
then the correlative pivot $\hat{u}(c)$ of $c$
lies on the quadric $Q$. 

{\bf (ii)}
Let $c \in \hat{F} - \hat{\Omega}$, and let 
$\hat{u} = \hat{u}(c)$ be the correlative pivot of $c$.
Then the set 
$$
q(c) =  \{ x \in {\bf P}^5:  Q = Q_x  \subset H_c \}
$$
is a smooth conic in the lagrangian plane
${\bf P}^2_{\hat{u}}$ of $\hat{u}$,
which we call the {\it vertex conic} of $c$. 
}

\bigskip

{\bf Proof of (i)}. \
By the $Sp_3$-homogeneity we can assume that $x = e_1$, 
and then in coordinates $(u:X:Y:z) = (u:(x_{ij}):(y_{ij}):z)$ 
$$
Q_{e_1} = (uy_{11} = x_{22}x_{23} - x_{23}^2) \subset
{\bf P}^4_{e_1} ={\bf P}^4(u:x_{22}:x_{23}:x_{33}:y_{11}). 
$$


Any $H_c$ which contains $Q_x =  Q_{e_1}$
contains ${\bf P}^4_x = {\bf P}^4_{e_1}$.
Therefore
$$
H_c \supset Q_{e_1}
\Leftrightarrow
c \in {\bf P}^8_{e_1} := ({\bf P}^4_{e_1})^{\perp}
= (u = x_{12} = x_{13} = x_{23} = y_{11} = 0)
\subset \hat{\bf P}^{13}(u:X:Y:z).
$$

Now the straightforward check shows that:
\ {\bf 1}. 
If $c \in {\bf P}^8_{e_1}$
then $F(c) = 0$, see {\bf (2.5.a)}. Therefore $c \in \hat{F}$. 
\ {\bf 2}. 
If $c \in {\bf P}^8_{e_1}$ then all the partial derivatives
$F_{z}(c)$, $F_{x_{1i}}(c)$, $i = 1,2,3$
and $F_{y_{ij}}(c)$, $(i,j) \not= (1,1)$ vanish.
Therefore $\nabla_F(c) \in {\bf P}^4_{e_1}$;
and since $\nabla_F(c) = \hat{u}(c) \in \Sigma$
is the correlative pivot of $c$ (see {\bf (2.6)}), then
$\hat{u}(c) \in \Sigma \cap {\bf P}^4_{e_1} = Q_{e_1}$.
{\bf q.e.d.}

\

{\bf Proof of (ii)}. \
Let 
$\rho : Sp_3 \times {\bf P}^{13} \rightarrow {\bf P}^{13},
{\rho}_g: x \rightarrow {\rho}_g(x)$
be the projective action
of $Sp_3$ on ${\bf P}^{13}$ (see {\bf (2.2)}),
and let ${\rho}_1 : x \mapsto g(x)$
be the projectivized standard action $\wedge_1$
of $Sp_3 \subset GL_6$ on ${\bf P}^5$. 
By the transitivity of $\hat{\rho}$ on the orbit
$\hat{F} - \hat{\Omega}$, we can choose $c = (0:0:Y:0)$,
where $Y$ is any symmetric rank $3$ matrix
(see {\bf (2.3)-(2.4)}),
and let $Y$ be the matrix of the quadratic
form $2y_1y_2+y_3^2$.    
Then
$c = 2y_{12} + y_{33} = x_{416}+x_{256}+x_{453}$: mod.${\bf C}^*$,
and the pivot $u(c) = piv(c) = (0:0:0:1) = (x_{456})$. 
The correlative point of $c$ is
$\hat{c} = (e_{143} + e_{523} + e_{126}) \in F -  \Omega$
(see {\bf (2.5.c)}), and the pivot of $\hat{c}$ is    
$\hat{u} = piv(\hat{c}) =  (1:0:0:0) \in \Sigma$,  
with a lagrangian plane
${\bf P}^2_{\hat{u}} = {\bf P}^2(x_1:x_2:x_3) = {\bf P}(U_o)$.
Let
$$
PSt_{\hat{c}}
= \{ g  \in Sp_3 : \rho_g(\hat{c}) = \hat{c} \}/{\bf C^*}
\subset PSt_3
$$
be the projectivized stabilizer group of $\hat{c}$.


\newcommand{\lotimes}
{\begin{picture}(7,3)
\put(2,0){$\subset$}
\put(2.2,0){$\times$}
\end{picture}}

\newcommand{\rotimes}
{\begin{picture}(7,3)
\put(2,0){$\times$}
\put(2.2,0){$\supset$}
\end{picture}}


\bigskip

{\bf ($\ast$)} {\bf Lemma}.
{\sl
The group $PSt_{\hat{c}} \subset Aut \ H_c$; 
and the correlative pivot plane
${\bf P}^2_{\hat{u}}$
is invariant under the standard action $\rho_1$
of $PSt_{\hat{c}}$. Moreover $PSt_{\hat{c}}$
is a projectivized semi-direct product
${\bf P}(O_q(3) \lotimes G_a^5)$,
where $PO_q(3)$ is the projectivized 
orthogonal group of the quadric   
$q = q(c) = (2x_1x_2 + x_3^2 = 0) \subset  
{\bf P}^2_{\hat{u}} = {\bf P}^2(x_1:x_2:x_3)$,
and $G_a^5$ is the additive group of ${\bf C}^5$.  
}

\bigskip 

{\bf Proof}. See \S 9 (2) in \cite{K}, or (5.21) in \cite{KS}.  

\bigskip 

Let $Q = Q_x$ be a $3$-fold quadric in $H_c$. 
By {\bf (i)}, the vertex $x$ of $Q = Q_x$ lies
in ${\bf P}^2_{\hat{u}}$;
and by Lemma {\bf ($\ast$)}, for any $g \in PSt_{\hat{c}}$
the quadric ${\rho}_g(Q_x) =  Q_{g(x)}$ 
also lies in $H_c$.
Therefore $Q =  Q_x$ lies in $H_c$, together with
all the quadrics $Q_{g(x)}, g \in PSt_{\hat{c}}$.
Again by Lemma {\bf ($\ast$)},
the action $\rho_1: (g,x) \mapsto g(x)$
of $PSt_{\hat{c}}$ on ${\bf P}^2_{\hat{u}}$ 
coincides with the action of $PO_{q}(3)$,
$q = q(c)$. Therefore $PSt_{\hat{c}}$
has two orbits in
${\bf P}^2_{\hat{u}} = {\bf P}^2(x_1:x_2:x_3)$
-- the smooth conic
$q = q(c) = (2x_1x_2 + x_3^2=0)$,
and its complement ${\bf P}^2_{\hat{u}} - q(c)$. 

Therefore in order to prove {\bf (ii)} it is enough to
see that:
\ {\bf 1}. for some $x  \in q(c)$ the 5-fold $H_c$ contains
the quadric $Q_x$.
\ {\bf 2}. for some $x \in {\bf P}^2_{\hat{u}}$ the quadric
$Q_x$ does not lie  in $H_c$.

This is straightforward:
If $x = e_1 \in q(c)$ then the quadric
$Q_{e_1}$ lies in $H_c = \Sigma \cap (2y_{12}+y_{33}=0)$, 
while for the point $x = e_3 \in {\bf P}^2_{\hat{u}} - q(c)$
the quadric $Q_{e_3}$ does not lie in $H_c$. 
{\bf q.e.d.}

%

\bigskip

{\bf (3.4)} {\bf Conics on ${\bf \Sigma}$}. \ 
Call a {\it formal conic} on $\Sigma$ 
any plane ${\bf P}^2 \subset {\bf P}^{13} = Span \ \Sigma$
such that the intersection cycle
$q = {\bf P}^2 \cap \Sigma$ is of dimension $1$
and of degree $2$.
Let ${\cal F}({\Sigma}) \subset G(3,14)$
be the family of formal conics on $\Sigma$.
Since $\Sigma$ is an intersection
of quadrics (see {\bf (2.1)}$({\bf \ast})$), 
and $\Sigma$ does not contain planes
(see e.g. Lemma 2.5.1 in \cite{IR}),
then the definition is correct.
The rank of the formal conic ${\bf P}^2$ is the
rank of the conic
$q = {\bf P}^2 \cap \Sigma \subset {\bf P}^2$. 
If $rank(q)$ = 2 or 3, then one can identify
the  conic $q \subset \Sigma$ 
and the formal conic ${\bf P}^2 = Span(q)$. 

\bigskip

{\bf (3.5)} {\bf Lemma}. \
{\sl
If $q \subset \Sigma$ is a conic of rank $\ge 2$,
then there exists a unique point $x = x(q) \in {\bf P}^5$,
such that $q \subset Q_x$.
}

\

We call the point $x = x(q)$ the {\it vertex} of
the $rank \ge 2$ conic $q$. 

\bigskip

{\bf Proof}. \ 
Since $q \subset \Sigma \subset G(3,6)$ is a conic
then the union $Q(q) := \cup_{u \in q} \ {\bf P}^2_u$
is either a ${\bf P}^3$ or a 3-fold quadric.

Suppose first that $Q(q) = {\bf P}^3$,
i.e. all the ${\bf P}^2_u, u \in q$ lie in ${\bf P}^3$;  
and let $u,v \in q$, $u \not=v$.
Since the lagrangian planes ${\bf P}^2_u$ and ${\bf P}^2_v$
both lie in ${\bf P}^3$, then they will intersect
each other along an isotropic line $L$.
By {\bf (3.1)}{\bf (i)}, the line $l = l_L$ with axis $L$
contains the points $u$ and $v$, therefore
$l = Span(u,v)$ and ${\bf P}^3_L = {\bf P}^3$.
But by {\bf (3.1)}{\bf (ii)} the only lagrangian planes 
in ${\bf P}^3_L$ are ${\bf P}^2_w, w \in l$. 
Therefore $q \subset l$ (as sets), which is only possible
if $rank(q) = 1$ -- contradiction.

Therefore $Q(q) \subset {\bf P}^4$ is a quadric,
and since on $Q(q)$ lie planes, $Q(q)$ is singular.

If $q = l+m$ is of rank $2$ then
$Q(q)$ is the union of the spaces ${\bf P}^3_L$
and ${\bf P}^3_M$ of the lines $l$ and $m$.
Let $L \subset {\bf P}^3_L$ and $M \subset {\bf P}^3_M$
be the axes of $l$ and $m$. Since $l \cap m = u$
is a point, then $Span(L \cup M) = {\bf P}^2_u$.
Therefore $L \cap M = x =  x(q)$ is a point, and
$q = l+m \subset Q_x$.

If $rank(q) = 3$ then $Q(q)$ is irreducible,
and since on $Q(q)$ lie planes then $Q(q)$
is singular and all the ${\bf P}^2_u, u \in q$
pass through the subspace ${\bf P}^k = Sing \ Q(q)$,
$k =  0,1$. If $k =  1$ then ${\bf P}^k = L$ is
a line. The line $L$ is isotropic since it lies
on lagrangian planes, and let $l  = l_L$ be the
line with axis $L$. But then 
${\bf P}^2_u \supset L, \forall u \in q$, 
together with {\bf (3.1)}{\bf (i)},  
will imply $q \subset l$ -- contradiction. 
Therefore $k = 0$, $x = x(q) := {\bf P}^k = {\bf P}^0$ 
is a point, and $q \subset Q_x$. 
{\bf q.e.d.}

\bigskip

{\bf (3.6)} {\bf Lemma}. \ 
{\sl
Let $l \subset \Sigma$ be a line.
Then the set ${\bf P}_l \subset G(3,14)$
of these planes ${\bf P}^2 \subset {\bf P}^{13}$
for which the intersection cycle
${\bf P}^2 \cap \Sigma = 2.l$,
is a line in $G(3,14)$. 
}

\bigskip

{\bf Proof}. \ 
Let ${\bf P}^2 \cap \Sigma  = 2.l$, and let $(u,v)$
be any pair of non-coincident points of $l$. 
Since ${\bf P}^2 \cap \Sigma  = 2.l$
then ${\bf P}^2 \subset {\bf P}^6_u \cap {\bf P}^6_v$, 
where ${\bf P}^6_u$ and ${\bf P}^6_v$
are the tangent projective spaces to $\Sigma$
at $u$ and $v$, see {\bf (2.3)-(2.4)}.
As in the proof of {\bf (3.1)}, we may assume that
$l = l_L = Span(u,v) = Span(e_{123}, e_{423})$
is the line with axis $L = Span(e_2,e_3)$.
Therefore ${\bf P}^3_l = {\bf P}^6_u \cap {\bf P}^6_v$
= $Span(e_{123},e_{423},e_{143}+e_{523},e_{124}+e_{623})$
is a projective 3-space containing the line
$l = Span(e_{123},e_{423})$.
Now it is clear that
${\bf P}^3_l = \cap_{u \in l} \ {\bf P}^6_u$,
and the set ${\bf P}_l$ coincides with the Schubert line
${\sigma}_{11,11,10}(l,{\bf P}^3_l)$ =
$\{ {\bf P}^2 \subset {\bf P}^{13} :
l \subset {\bf P}^2 \subset {\bf P}^3_l \}$
$\subset G(3,14)$.
{\bf q.e.d.}

\bigskip

{\bf (3.7)} {\bf Corollary}. \ 
{\sl
Let
${\cal F}_k({\Sigma}) \subset {\cal F}({\Sigma})  \subset G(3,14)$
be the loci of formal conics on $\Sigma$ of rank $\le k$,
$k = 1,2$. 
Then $dim\ {\cal F}({\Sigma})=11$,
$dim\ {\cal F}_2({\Sigma})=10$, 
and \ $dim \ {\cal F}_1({\Sigma})=8$.  
}

\bigskip

{\bf Proof}. \ 
By {\bf (3.4)} any conic $q$ of rank $>1$
can be identified with its formal conic
${\bf P}^2(q) =  Span(q)$.
By {\bf (3.5)} any such $q$ is a linear section of
a unique quadric $Q_{x(q)}$ with the plane
${\bf P}^2(q) \subset{\bf P}^4_{x(q)} = Span \ Q_{x(q)}$. 
Therefore $dim \ {\cal F}({\Sigma})$ 
= $dim \ {\bf P}^5 + dim \ G(3,5) = 11$, 
and $dim \ {\cal F}_2({\Sigma}) = 10$.
At the end, {\bf (3.1)} and {\bf (3.6)} imply:
$dim \ {\cal F}_1({\Sigma}) = dim \ LG_2 + 1 = 8$,
see {\bf (2.1)}.
{\bf q.e.d.}

\bigskip
\bigskip 
\bigskip

\centerline{\large \bf \S 4.
 The Fano family ${\bf {\cal F}(X)}$ of conics on ${\bf X}$}

\bigskip
\medskip

{\bf (4.1)}
{\bf The dual plane quartic ${\bf C_X}$ of ${\bf X}$}. \
By \cite{M1}, \cite{M2} any smooth prime Fano 3-fold
$X = X_{16}$ of degree 16 is a linear section of the
$Sp_3$-grassmannian $\Sigma \subset {\bf P}^{13}$
by a codimension 3 subspace
${\bf P}^{10} \subset {\bf P}^{13}$.
Moreover two prime Fano 3-folds
$X_{16}' = \Sigma \cap {\bf P}^{10}_1$
and $X_{16}'' = \Sigma \cap {\bf P}^{10}_2$
are projectively equivalent \ {\it iff} \ 
${\bf P}^{10}_1$ and ${\bf P}^{10}_2$ are conjugate
under the action $\rho$ of $Sp_3$ in ${\bf P}^{13}$.

\

From now on we shall consider only the situation
when $X$ is {\it general}.

\

Fix a general prime Fano 3-fold
$X = X_{16} = \Sigma \cap {\bf P}^{10}_X$; 
and let
${\bf P}^2_X = {{\bf P}^{10}_X}^{\perp} \subset \hat{\bf P}^{13}$
be the plane of linear equations of
${\bf P}^{10}_X \subset {\bf P}^{13}$.
Since $X$ is general then the plane
${\bf P}^2_X = {{\bf P}^{10}_X}^{\perp} \subset \hat{\bf P}^{13}$
of linear equations of ${\bf P}^{10}_X \subset {\bf P}^{13}$
intersects the invariant quartic hypersurface $\hat{F}$ 
along a smooth plane quartic curve
$$
C_X = \hat{F} \cap {\bf P}^2_X.  
$$
By the preceding, $C_X$ is uniquely defined by the
Fano 3-fold $X$; and we call the curve $C_X \subset {\bf P}^2_X$
{\it the dual plane quartic} of $X$.

Since $X$ is general then the plane ${\bf P}^2_X$ does not
intersect the codimension $3$ orbit $\hat{\Omega}$.
Therefore any point $c \in C_X$
has a uniquely defined double pivot $u = u(c) \in \hat{\Sigma}$. 
By {\bf (2.7)}{\bf (i)}-{\bf (ii)},  
the hyperplane sections $H_c$ and $H_u$ of $\Sigma$ both
contain the Veronese cone $K_{\hat{u}} = {\bf P}^6_{\hat{u}} \cap \Sigma$.
Therefore any point $c \in C_X$ defines uniquely a $4$-fold linear
section
$$
W_c = H_c \cap H_{u(c)} \subset \Sigma; 
$$
and since $K_{\hat{u}(c)} =  Sing \ H_u$ (see {\bf (2.7)}{\bf (ii)})
then $W_c$ is singular along $K_{\hat{u}(c)}$.

\bigskip

{\bf (4.2)} {\bf Lemma}. \ 
{\sl
Let $c \in \hat{F} - \hat{\Omega}$, and let
$q(c) \subset {\bf P}^2_{\hat{u}(c)}$ be the
vertex conic of $c$, see {\bf (3.3)}. Then
$W_c = \cup_{x \in q(c)} \ Q_x$.
}

\bigskip

{\bf Proof}. \
By {\bf (2.7)} and {\bf (3.2)}
$$
H_{u(c)} = {\sigma}_{100}({\bf P}^2_{\hat{u}(c)}) \cap {\bf P}^{13}
= \{ w \in \Sigma :
   {\bf P}^2_w \cap {\bf P}^2_{\hat{u}(c)} \not = \emptyset \}
= \cup_{x \in \ {\bf P}^2_{\hat{u}(c)}} \ Q_x .
$$
Moreover, by {\bf (3.3)} any quadric
$Q_x, x  \in q(c) \subset {\bf P}^2_{\hat{u}(c)}$ 
lies in $H_c$. {\bf q.e.d.}

\bigskip

{\bf (4.3)}
{\bf The singular hyperplane sections ${\bf S_c \subset X}$}. \
For the general $X$, the dual curve $C_X$ does not intersect
$\hat{\Omega}$, hence the pivot map
$piv: C_X \rightarrow \hat{\Sigma}, \ c \mapsto u(c) = piv(c)$
is regular.
The regular map $piv :C_X \rightarrow piv(C_X)$
is an isomorphism. Indeed if $c_1,c_2 \in  C_X$
are two points such that  $u(c_1) = u(c_2) = u$,  
then both $c_1,c_2$ will lie in the $6$-space
${\bf P}^6_u$ (see {\bf (2.5.b)} and {\bf (2.3)(ii)}); 
and since ${\bf P}^6_u  \subset \hat{F}$ (see {\bf (2.3)(ii)})
then the line $L = Span(c_1,c_2) \subset {\bf P}^2_X$
will lie in $\hat{F}$. But then the line $L$
will be a component of $C_X = \hat{F} \cap {\bf P}^2_X$,
which contradicts the general choice of $X$. 
%
%
Since the plane ${\bf P}^2_X = Span \ C_X$ does not intersect
$\hat{\Sigma} \supset piv(C_X)$ then for any $c \in C_X$ the set 
$$                                       
S_c =  W_c \cap {\bf P}^{10}_X = X \cap H_{u(c)} \subset X,
$$
where $u(c) = piv(c)$, is a hyperplane section of $X$.
Since different points $c_1,c_2 \in  C$ have different
pivots $u(c_1),u(c_2) \in piv(C_X)$ then
$S_{c_1} \not= S_{c_2}$ for $c_1 \not= c_2$. 
By construction any surface $S_c,  c \in C_X$
is a complete intersection of two hyperplanes in $W_c$,
say $S_c = W_c \cap H_1 \cap H_2$.  
Since the $5$-fold $W_c$ is singular along the Veronese cone
$K_{\hat{u}(c)}$ (see {\bf (4.1)}), 
then the surface $S_c =  W_c \cap H_1 \cap H_2$
is singular along the intersection cycle
$$
C_c = K_{\hat{u}(c)} = K_{\hat{u}(c)} \cap H_1 \cap H_2 
= K_{\hat{u}(c)} \cap {\bf P}^{10}_X.
$$

Assume that $\hat{u}(c) \in C_c$.
Then $\hat{u}(c) \in H_1 \cap H_2$; 
and since $X = H_c \cap H_1 \cap H_2$
and $H_c$ is singular at $\hat{u}(c)$ (see {\bf (2.7)}{\bf (i)}) 
then $X$ will be singular at $\hat{u}(c)$,
which contradicts the general choice of $X$.
Therefore $C_c$ is a codimension 2 linear section
of the Veronese cone $K_{\hat{u}(c)}$, which
does not contain the vertex $\hat{u}(c)$ of $K_{\hat{u}(c)}$,
i.e. $C_c$ is a rational normal quartic curve
(including the case when $C_c$ is a union
of two conics intersecting each other at a point).

\bigskip

{\bf (4.4)} {\bf Lemma}. \ 
{\sl
Let $c \in  C_X$ and let  $q(c) \subset {\bf P}^2_{\hat{u}(c)}$
be the vertex conic of $c$. Then for any $x \in q(c)$
the intersection $q_x  = Q_x \cap S_c$ is a conic.
If $c \in C_X$ is general then for the general
$x \in q(c)$ the conic $q_x$ is a bisecant to the rational
normal quartic $C_c \subset Sing \ S_c$.
Moreover for any $c \in C_X$ the surface
$S_c = \cup_{x \in q(c)} \ q_x$. 
}

\bigskip

{\bf Proof}. \
Let $c \in C_X$, etc. be as above.
We shall see first that for any $x \in q(c)$
the quadric $Q_x \subset W_c$ (see {\bf (4.2)}) intersects
the Veronese cone $K_{\hat{u}(c)}$
along a cone $K^x_{\hat{u}(c)} \subset K_{\hat{u}(c)}$
over a conic on the Veronese surface.
Indeed the intersection
$K^x_{\hat{u}(c)}$ = $Q_x \cap K_{\hat{u}(c)}$
=    $\{ w \in \Sigma : x \in {\bf P}^2_w \ \& \
     dim \ ({\bf P}^2_w \cap {\bf P}^2_{\hat{u}}) \ge 1 \}$
     coincides with the union of lines 
%
%
$\cup \ \{ l_L : x \in L \subset {\bf P}^2_{\hat{u}} \}$
\ on $\Sigma$, see {\bf (3.1)}. 
%
%
Any such $l_L$ is a ruling line of the cone $K_{\hat{u}(c)}$.
Therefore $K^x_{\hat{u}(c)} \subset K_{\hat{u}(c)}$
is a cone; and since the base set
$\{ L \subset {\bf P}^2_{\hat{u}}: x  \in L \}$
is a line in $\hat{\bf P}^2_{\hat{u}}$ then the base
of $K^x_{\hat{u}(c)}$ is a conic on the Veronese surface.
%
%
By {\bf (4.2)}, for any $x \in q(c)$ the quadric $Q_x$ lies in $W_c$; 
and since
$X = W_c \cap {\bf P}^{10}_X = W_c \cap  H_1 \cap  H_2$ 
then the cycle 
$q_x := Q_x \cap {\bf P}^{10}_X$
= $Q_x \cap H_1 \cap H_2  \subset W_c \cap {\bf  P}^{10}_X = X$
will be a quadric of dimension at least $1$.
Since $Pic(X) = {\bf Z}.H$, where $H$ is the hyperplane section,
then the prime Fano $3$-fold $X  = X_{16}$ can't contain
quadric surfaces.

Therefore for any $c \in C_X$, and any $x \in q(c)$,
the cycle $q_x$ is a conic; and 
$q_x$ intersects the rational quartic
$C_c \subset Sing \ S_c$ at the cycle
$$
z_x  =  q_x \cap C_c = K^x_{\hat{u}(c)} \cap H_1 \cap H_2.
$$
By construction, the $4$-fold $W_c$
and the $rank \ 3$ quadratic cone
$K^x_{\hat{u}(c)}$ = $Q_x \cap K_{\hat{u}(c)}$
depend only on the choice of the point
$c \in \hat{F} - \hat{\Omega}$
and the point $x$ on its vertex conic $q(c)$.
Therefore for the general codimension
$2$ linear section $X = W_c \cap H_1 \cap H_2$, 
the cycle $z_x = K^x_{\hat{u}(c)} \cap H_1 \cap H_2$
will be a general section
of $K^x_{\hat{u}(c)}$ with two hyperplanes,
and then $q_x \subset X$ will be a smooth conic
on $S_c$ bisecant to the rational normal quartic
$C_c$.
For arbitrary $x \in q(c), c \in C_X$,
the cycle $z_x$ is a linear section, of dimension $\le 1$,
of the quadratic cone $K^x_{\hat{u}(c)}$. 
In the extremal case when $dim \ z_x = 1$, 
the curve $C_c$ must be a union of two conics
on the Veronese surface such that $z_x = q_x$
is one of these two conics.
This can be possible for at most two conics
from the pencil $\{  q_x : x \in q(c) \}$.
In all the rest possible cases $z_x$ will be 
a zero-cycle of degree $2$.

At the end, the conics $q_x, x \in  q(c)$ sweep
the surface $S_c$ out, since
$S_c = W_c \cap H_1 \cap  H_2$
= $\cup_{x \in q(c)} \  Q_x  \cap H_1 \cap H_2
= \cup_{x \in q(c)} \ q_x$, see {\bf (4.2)}. 
{\bf q.e.d.}

%

\bigskip

{\bf (4.5)}
{\bf The Fano surface ${\bf {\cal F}(X)}$
and the vertex surface ${\bf S_X}$}. \ 
For the general $X$ the family ${\cal F}(X)$
of conics on $X$ is a union of $2$-dimensional
components ${\cal F}^{1}, ..., {\cal F}^{n}$,
see \S 4 in \cite{I1} or \S 4.2 in \cite{IP}.
We shall see that $n = 1$, and ${\cal F}(X) = {\cal F}^{1}$
is a ruled surface over the dual plane quartic $C_X$ of $X$. 

\bigskip

{\bf (4.5.a)}  
{\bf Lemma}. \
{\sl
The general $X = X_{16}$ does not contain conics of rank $1$.
}

\bigskip

{\bf Proof}. \
For the general $X$, the codimension $9$ Schubert cycle
${\sigma}_{333}({\bf P}^{10}_X) = G(2,{\bf P}^{10}_X)$,
in $G(3,14) = G(2:{\bf P}^{13})$, 
does not intersect the $8$-fold
${\cal F}_1(\Sigma) \subset G(3,14)$ of formal conics
of rank $1$ on $\Sigma$, see {\bf (3.7)}.
{\bf q.e.d.}

\bigskip

{\bf (4.5.b)} {\bf The vertex map}. \  
By {\bf (3.5)} and {\bf (4.5.a)}, for the
general $X$ the map
$$
ver:  {\cal F}(X) \rightarrow {\bf P}^5,
\  q \mapsto \ ver(q) = \mbox{ the vertex } \ x(q) \ \mbox{ of } \ q
$$
is regular; and we call this map the
{\it vertex map} of ${\cal F}(X)$. 

\bigskip

{\bf (4.5.c)} 
{\bf Lemma}. \
{\sl The vertex map $ver$ sends the family
${\cal F}(X)$ isomorphically onto 
the irreducible closed set 
$$
S_X = \cup_{c \in C_X} \ q(c) \subset {\bf P}^5,
$$
swept out by the smooth vertex conics
$q(c) \subset {\bf P}^2_{\hat{u}(c)}$
of the points $c \in C_X$, see {\bf (3.3)}.
}

\bigskip 

{\bf Proof}. \
Let $x \in S_X$. Then $\exists c \in C_X$
such that $x \in q(c)$. By {\bf (4.4)}
the cycle $q_x = Q_x \cap S_c$ is a conic on
$X$ with vertex $x$.
Therefore $S_X \subset ver({\cal F}(X))$.

Next, we shall see that
$ver({\cal F}(X)) \supset S_X$. 
For this, let $q \subset X$ be a conic,
and let $x =  ver(q)$.
Since $q \subset Q_x$ is a codimension $2$
linear section of the 3-fold quadric $Q_x$
(see {\bf (3.5)}), 
and since any hyperplane section
$H_{c}$, $c \in {\bf P}^2_X$ contains $q \subset X$,
then the plane ${\bf P}^2_X$ contains a
{\it unique} point $c =  c(x)$ such 
that $Q_x \subset H_c$.
Since $H_c$ contains a 3-dimensional quadric
(the quadric $Q_x$) then, by {\bf (3.3)}{\bf (i)},
the point $c$ must lie in $\hat{F}$. 
Therefore $c \in {\bf P}^2_X \cap \hat{F} = C_X$.
Moreover, by {\bf (3.3)}{\bf (ii)},
the vertex $x = x(q) = ver(q)$ must lie in the vertex
conic $q(c) \subset {\bf P}^2_{\hat{u}(c)}$ of $c$;
and since $q(c) \subset S_X$ then
$x = ver(q) \in q(c)$ lies in $S_X$.  
Therefore $ver({\cal F}(X)) \supset S_X$. 
Moreover, by the proof of {\bf (4.4)},
the conic $q_x$ is the unique conic
on $X$ with vertex $x \in q(c)$.
Therefore the regular and surjective
vertex map $ver: {\cal F}(X) \rightarrow S_X$
is injective.
{\bf q.e.d.}

\bigskip

{\bf (4.5.d)} {\bf The vertex surface ${\bf S_X}$}. \ 
By {\bf (4.5.a-c)}, the vertex map $ver$
sends ${\cal F}(X)$ isomorphically
onto the irreducible closed set 
$$
S_X  = \cup_{c \in C_X} \ q(c).
$$
Therefore ${\cal F}(X)$ is irreducible, 
and $S_X \subset {\bf P}^5$ is a surface 
which we call the {\it vertex surface} of $X$.

\bigskip

Lemma {\bf (4.4)} and {\bf (4.5.a-d)} imply the following 

\bigskip

{\bf (4.6)} {\bf Theorem}. \
{\sl
Let $X  =  X_{16}$ be general. Then the Fano family ${\cal F}(X)$
of conics on $X$ is a ruled surface
$p: {\cal F}(X) \rightarrow C_X$  
over the dual plane quartic curve
$C_X$ of $X$. For the point $c \in C_X$, the fiber
$f_c = p^{-1}(c) \subset {\cal F}(X)$
coincides with the pencil $\{ q_x, x \in q(c) \}$
of conics $q_x \subset X$ defined in {\bf (4.4)}. 
The vertex map $ver: {\cal F}(X) \rightarrow {\bf P}^5$
sends ${\cal F}(X)$ isomorphically onto the vertex surface
$S_X = \cup_{c\in C_X} \ q(c) \subset {\bf P}^5$
of $X$. For any point $c \in C_X$, the isomorphism 
$$
ver: {\cal F}(X) \rightarrow S_X \subset {\bf P}^5
$$
sends the fiber $f_c \subset {\cal F}(X)$
onto the correlative pivot conic
$q(c) \subset {\bf P}^2_{\hat{u}(c)}$ of $c$.
}

\bigskip

{\bf Proof}. \
By {\bf (4.4)} and {\bf (4.5.c)}, it only rests
to see that the vertex surface $S_X \cong {\cal F}(X)$
is ruled over $C_X$,
with fibers -- the conics $q(c), c  \in C_X$.
For this it is enough to see that if
$c_1,c_2 \in C_X$, $c_1 \not=  c_2$
then the conics $q(c_1)$ and $q(c_2)$
do not intersect each other.

Suppose that $q(c_1) \cap q(c_2) \not= \emptyset$.
Let $x \in q(c_1) \cap q(c_2)$,
and let $H_{c_1}$ and $H_{c_2}$ be the hyperplane sections
of $\Sigma$ defined by $c_1$ and $c_2$. 
Since $x \in S_X$ then, by the proof of {\bf (4.4)}, 
there exists a unique conic $q = q_x \subset X$
with vertex $x$. Let also $Q_x \subset \Sigma$ be
the $3$-fold quadric with vertex $x$; in particular 
$q_x = Q_x \cap {\bf P}^{10}_X$, ibid. 
By {\bf (3.3)(ii)}, both $H_{c_1}$ and $H_{c_2}$
must contain $Q_x$ since the vertex $x$ of $Q_x$
lies in their pivot conics $q(c_1)$ and $q(c_2)$. 
But by the proof of {\bf (4.5.c)}, there exists
a unique point $c = c(x) \in C_X$ such that
$Q_x \subset H_c$. Therefore $c_1 = c_2 = c$.
{\bf q.e.d.}

\newpage

\centerline{\large \bf \S 5.
 Prime Fano 3-folds of degree 16 and stable ruled surfaces}

\smallskip

\centerline{\large \bf of invariant 3 over plane quartics}

\bigskip
\medskip

{\bf (5.1)}
{\bf Lemma}:
{\bf the double projection from a line ${\bf l \subset X_{16}}$}
(V. Iskovskikh).
{\sl

{\bf (1)} \
Let $X = X_{16} \subset {\bf P}^{10}$
be a smooth prime Fano $3$-fold of degree $16$,
and let $l \in X$ be a line. 
Then the double projection ${\pi}_{2.l}$
from the line $l$, given by the non-complete linear system
$\mid {\cal O}_X(1 - 2.l) \mid$,
defines a birational isomorphism
${\pi} = {\pi}_{2.l}:  X \rightarrow {\bf P}^3_l$,
where ${\bf P}^3_l$ is the $3$-dimensional
projective space. 

There exists a smooth curve $C^l \subset {\bf P}^3_l$
of genus $3$ and of degree $7$, 
which lies on a unique cubic surface $S = S_3$, such that
the inverse birational map ${\varphi}: {\bf P}^3_l \rightarrow X$
is defined by the non-complete linear system
$\mid {\cal O}_{{\bf P}^3_l}(7 - 2{C^l})\mid$, and:

{\bf (i)} The one-dimensional family ${\cal Q}_l$
of conics $q \subset X$ which intersect $l$ 
sweeps out the unique effective divisor
$Q = Q_l$ from the linear system $\mid {\cal O}_X(3 - 7.l) \mid$.

{\bf (ii)} The double projection ${\pi} = {\pi}_{2.l}$
can be represented as a product
$\pi = {\tau} \circ {\rho} \circ {\sigma}^{-1}$,
where ${\sigma}:X' \rightarrow X$ is the blowup of
$l \in X$, ${\rho}:X \rightarrow X^+$ is
a flop over the projection ${\pi}_l:X \rightarrow X''$ from $l$,
and ${\tau}:X^+ \rightarrow {\bf P}^3_l$
is a blow-down of the proper image
$Q^+ \subset X^+$ of $Q$ onto the curve 
${C^l} \subset {\bf P}^3_l$. 
The extremal curves $q^+ \subset X^+$ contracted by $\tau$ 
are the strict transforms of the conics
$q \in {\cal Q}_l$.

{\bf (iii)} The unique cubic surface
$S = S_l \subset {\bf P}^3_l$ through the curve $C^l$, 
is swept out by the one-dimensional
family ${\cal S}_l$ of conics $s \subset {\bf P}^3_l$,
intersecting $C^l$ at a $0$-cycle of degree $7$. 
The proper transform $S' \subset X'$
of $S$ coincides with the exceptional divisor
${\sigma}^{-1}(l) \subset X'$ of $\sigma$.
The strict transforms $s' \subset X'$ of the conics
$s \in {\cal S}_l$ are the extremal curves of $\sigma$.

{\bf (iv)} There exists a non-negative integer
$e  = e(l) \le 5$ (and if $X$ is general and
$l \subset X$ is general then $e(l) = 5$)
such that the flop $\rho : X' \rightarrow X^+$
transforms the proper $\sigma$-preimages
$l'_1,...,l'_e \subset X'$ of the $e$ lines
$l_1,...,l_e \subset X$ which intersect $l$
to the proper $\tau$-preimages $l^+_i$, $i = 1,..,e$
of $e$ lines $L_1,...,L_e \subset {\bf P}^3_l$ of $C^l$.
If $X$ is general and $l \subset X$ is general,
then the lines $L_1,...,L_5$ are the five
$4$-secant lines to $C^l$. 

{\bf (2)} \ Let $C \subset {\bf P}^3$  be a smooth curve 
of genus $3$ and of degree $7$, which lies on a unique
cubic surface $S = S_3$. Then there exists a smooth prime
Fano $3$-fold $X$ of degree $16$ and a line $l \subset X$,
such that {\bf (1)} takes place for $l \subset X$
and $C^l = C \subset {\bf P}^3_l =  {\bf P}^3$.
}

\bigskip 

{\bf Proof.} See \cite{I2}. 

\bigskip

{\bf (5.1.a)} {\bf Remark}. \
The general curve $C \subset {\bf P}^3$
of genus $3$ and degree $7$ evidently
lies on a unique and smooth cubic surface $S = S_3$.
Therefore $S = S_l$ and $C = C^l$ for some $l \subset X$
as in {\bf (5.1)}{\bf (2)}. 
One can represent (non-uniquely) the smooth cubic
$S_3 \subset {\bf P}^3$ as the blowup 
$S_3 = \tilde{\bf P}^2_{z_0,...,z_5}$ of
${\bf P}^2$ at $6$ points $z_0, ..., z_5$, and
the curve $C \subset S_3$ -- as the proper preimage
of an element of the system
$|{\cal O}_{{\bf P}^2}(4-z_1-...-z_5)|$.
Then the five $4$-secant lines of $C$ are
the proper preimages $L_i \subset S_3$ of the
conics $C_i \subset {\bf P}^2$
passing through the points $z_0,...,\hat{z_i},...,z_5$,
$i = 1,...,5$. 

Notice that the line $L = L_0$ = the blowup of the
point $z_0$ is the unique line on $S_3$ which does
not intersect $C$, and $L_1,...,L_5$ are the
five lines on $S_3$ which intersect $L$. 

\bigskip

{\bf (5.1.b)}  
{\bf Corollary}. \
{\sl The family $\Gamma(X)$ of lines on the general prime
Fano 3-fold $X$ is a smooth irreducible curve of genus $17$.
}

\bigskip

{\bf Proof}. \
By Th. 4.2.7 in \cite{IP}, for the general $X = X_{16}$
the $1$-dimensional family $\Gamma(X)$ of lines on $X$ is smooth.  
It rests to see that  $\Gamma(X)$ is irreducible and of genus $17$.

Let ${\pi} = {\pi}_{2.l}: X \rightarrow {\bf P}^3_l = {\bf P}^3$
be the double projection of $X$ from the general line $l \subset X$
as in {\bf (5.1)}.
If $m \subset X$ is a line which  does not intersect $l$
then its proper image ${\pi}_*(m)$ is evidently a line
in ${\bf P}^3$; and let ${\pi}_*(m)$ intersects $C^l$ at $k$ points.
Since by {\bf (5.1)}{\bf (1)}, $\varphi = \pi^{-1}$
is defined by the system
$\mid {\cal O}_{{\bf P}^3}(7 - 2C^l) \mid$,
then the proper image ${\varphi}_*({\pi}_*(m)) = m$ 
will be a curve of degree = $7.deg\ m(l) - 2k = 7 - 2k$; 
and since $m$ is a line then $k = 3$.
The same argument in the opposite direction gives
that any purely $3$-secant line to $C^l$ is the proper
$\pi$-image of a line $m \subset X$ which does not
intersect $l$. 

If $m = l_i, i = 1,..,5$ is one of the five $4$-secant lines
of $l$ then, by {\bf (5.1)}{\bf (iv)}, its proper image
${\pi}_*(m)$ in ${\bf  P}^3$ is
one of the five $4$-secant lines $L_i$ to $C^l$.  
These $L_i$ evidently are 4-tuple singularities
of the curve
$Sec_3(C^l)$ of the $3$-secant lines of $C^l$.
Therefore the map
$$
\pi_*:  \Gamma(X) \rightarrow Sec_3(C^l),  
$$
is the normalization of $Sec_3(C^l)$ 
at its five 4-tuple points $L_1,...,L_5$.
By {\bf (5.1)}{\bf (2)}, for the general $X$ the curve $C^l$
is a general curve of genus $3$ and degree $7$ in ${\bf P}^3$.
Therefore, by Proposition 2.4 and Theorem 3.6 in \cite{GP},
$\Gamma(X)$ is an irreducible smooth curve of genus $17$.
{\bf q.e.d.}

\bigskip       

{\bf (5.2) }
{\bf Deformations and minimal sections of ruled surfaces}
(see \cite{LN}, \cite{S}). \ 
The ruled surface $S = {\bf P}(F) \rightarrow C$ over the
smooth curve $C$ is stable if the rank $2$ vector bundle
$F \rightarrow C$ is stable. Since the stability of $F$
does not depend on the twist $F \otimes L$ by a line
bundle $L$, the definition is correct.

For the section $C_o \subset S$, denote by
$\varepsilon(C_o) =  C_o^2 = {\cal O}_{C_o}(C_o) \in Pic(C_o)
= Pic(C)$ its self-intersection divisor, and let  
$e(C_o) = deg \ \varepsilon(C_o) = deg \ C_o^2$.   
The integer {\it invariant}
$$
e(S) \ := min \ \{ e(C_o) : C_o \mbox{ is a section of } S \}
$$
is always $\le g = g(C)$. For a fixed $g$, the ruled surfaces $S$
form two deformation classes --
{\it even}:  with $e(S) \equiv g$ (mod.2),
and {\it odd}: with $e \equiv g-1$ (mod.2).
The versal deformations in the even (respectively odd) class
are ruled surfaces $S$ with $e(S) = g$
(respectively with $e(S) = g-1$).
For $g \ge 2$ the general $S$ of any of the
two versal types -- the even and the odd -- is stable;
and if ${\cal S}[g]$ and ${\cal S}[g-1]$ denote
respectively their moduli spaces, then  
$dim \ {\cal S}[g] = dim \ {\cal S}[g-1] = 6g-6$.
The same as above takes place if the smooth
base curve $C$ of genus $g \ge 2$ is fixed, with the only
difference that then the moduli spaces ${\cal S}_C[g]$
and ${\cal S}_C[g-1]$ of versal even and odd ruled surfaces
over $C$ both have dimension $3g-3$, see \cite{S}.
In particular, if $g(C) = 3$ then the ruled surfaces
of invariant $3$ over $C$ are exactly the elements 
of the even versal class, and their moduli space 
${\cal S}_C[3]$ has dimension $6$.

Let $S = {\bf P}(F) \rightarrow C$ be a stable ruled surface
with invariant $e(S) = e$, and let $C_o \subset S$
be the section defined by the extension

\medskip

{\bf ($\ast$)} \hspace{4.0cm}
$0\rightarrow{\xi}\rightarrow F\rightarrow{\eta}\rightarrow 0$ 

\medskip

\noindent for $\eta, \xi \in Pic(C)$, 
see Proposition 2.6 in Ch. 5 of \cite{H}. 
In particular
$\varepsilon = \varepsilon(C_o) = \eta - \xi$,
with operations in $Pic(C)$ written
additively.

The {\it minimal sections}
of $S$ are these sections $C_o$ of $S$ for which
$e(C_o) = deg \ \varepsilon(C_o) = e(S)$. 
Equivalently, the section $C_o \subset S = {\bf P}(F)$
is minimal \ {\it iff} \ the subbundle
${\xi} \subset F$ is of the maximum possible degree,
or a {\it maximal} {\it subbundle} of $F$,  
see Lemma 2.3 in \cite {LN}.
By Lemma 2.2 in \cite{LN} the minimal sections $C_o$ 
of $S$ are defined uniquely by their intersection
divisors $\varepsilon{(C_o)} = C_o^2$.
We shall denote by $Min(S)$ the family of minimal sections
of $S$, and by $C_{\varepsilon} \in Min(S)$ the minimal
section of self-intersection divisor
$\varepsilon(C_o) = \varepsilon \in Pic(C)$.

Let $g \ge 3$, and $g-1 \le e = e(S) \le g$.
For the fixed minimal section $C_{\varepsilon} \in Min(S)$,
$\deg \ \varepsilon = e = e(S)$, let
$[e_{\varepsilon}] \in {\bf P}^{g-2+e}_{\varepsilon}$
 = ${\bf P}((H^o(C, K_C + \varepsilon )^*)$
 = ${\bf P}(Ext^{1}({\eta},{\xi}))$
be the extension-class point defined by $C_{\varepsilon}$,
and let 
${\Phi}_{\varepsilon}: C \rightarrow {\bf P}^{g-2+e}_{\varepsilon}$
be the map defined by the complete linear system
$|K_C + \varepsilon|$.
By Proposition 2.4 in \cite{LN} there exists a
canonical bijection between the set
$Min(S) - \{ C_{\varepsilon} \}$ of minimal sections of $S$
different from $C_{\varepsilon}$ and the $e$-secant
spaces ${\bf P}^{e-1} \subset {\bf P}^{g-2+e}_{\varepsilon}$
to the curve ${\Phi}_{\varepsilon}(C)$ which pass through
$[e_{\varepsilon}]$. 

We shall see below that for $g = 3$, the general ruled surface
$S \in {\cal S}[3]$ is the same as the Fano surface
${\cal F}(X)$ of conics on a general prime Fano 3-fold $X = X_S$
of degree $16$ defined uniquely by $S$, and the family $Min(S)$
of minimal sections of $S$ is the same as the curve
$\Gamma(X_S)$ of lines on $X_S$.

\bigskip

{\bf (5.3)} {\bf Non-abelian Brill-Noether loci of type II}
(see \cite{M3}). \ 
Let $F \rightarrow C$ be a rank $2$ vector bundle over the smooth curve $C$,
let $K = K_C$ be the canonical bundle of $C$, and let $\nu$ be a nonnegative
integer such that ${\nu} \equiv {deg \ F} \ mod \ 2$.  The non-abelian
Brill-Noether locus of type II, associated to the triple $(C,K,F)$, is
the set of equivalence classes of bundles

\medskip

\centerline{$M_C(2,K:{n}F) = \{ E \rightarrow C:
        rank\ E = 2, det\ E = det\ F \otimes K, dim \ Hom(F,E) \ge {n} \}$.}

\medskip

It is proved in \S 6 of \cite{M3} that $M_C(2,K:{n}F)$ admits a
natural scheme structure as certain Pfaffian locus in the moduli space
$M_C(2, det\ F \otimes K_C)$ of stable rank $2$ vector bundles on $C$
of determinant $det\ F \otimes K_C$. 

\bigskip

{\bf (5.4)} \ {\bf Lemma} (S. Mukai). \
{\sl

{\bf (a)} \ Let $C \subset {\bf P}^3$ be a nonsingular plane quartic
curve, and let $F$ be a rank $2$ stable vector bundle on $C$
such that the ruled surface $S = {\bf P}(F) \rightarrow C$ 
has invariant $e = e(S) = 3$. 
Then the Brill-Noether locus of type II
$$
X = M_C(2,K:3F)
$$
is a nonsingular prime Fano $3$-fold of degree $16$. 

{\bf (b)} \ Let $C_{\varepsilon}$ be the minimal section of
$S$ of self-intersection divisor $\varepsilon$,
and let {\bf (5.2)}{\bf ($\ast$)}, with $\eta - \xi = \varepsilon$,
be the extension defined by $C_{\varepsilon}$. 
Then the locus
$$
l_{\varepsilon}
   = \{ E \in X = M_C(2,K:3F): h^o(E({-\eta})) \ge 2 \} \subset X
$$
is a line on the Fano $3$-fold $X = X_{16} \subset {\bf P}^{10}$.

{\bf (c)} \ 
In the notation of {\bf (5.2)}, the map 
${\Phi}_{\varepsilon}: C \rightarrow {\bf P}^4_{\varepsilon}$
is an embedding; and the projection 
$p_{[e_{\varepsilon}]}:{\bf P}^4_{\varepsilon}
        \rightarrow {\bf P}^3_{\varepsilon}$ \
from $[e_{\varepsilon}]$ sends ${\Phi}_{\varepsilon}(C)$
isomorphically onto a space curve $C^l_{\varepsilon}$ of
degree $7$.

{\bf (d)} \
The curve $C^l_{\varepsilon} \subset {\bf P}^3_{\varepsilon}$
is projectively equivalent to the curve 
$C^l \subset {\bf P}^3_l$ defining the inverse of the
double projection ${\pi}_{2.l}:X \rightarrow {\bf P}^3_l$,
see {\bf (5.1)}{\bf (1)}.
}

\bigskip

{\bf Proof}. See \S 9 in \cite{M3}.

\bigskip

{\bf (5.5)} \ 
{\bf The associate Fano 3-fold ${\bf X_S}$ of ${\bf S}$}. \
According to Remark 9.2 of \cite{M3}, the locus
$X_F = M_C(2,K:3F)$
does not depend on the twist $F \rightarrow F \otimes L$
by an invertible sheaf $L$ on $C$.
That is, the Fano $3$-fold
$X_S = X_F = M_C(2,K:3F)$ depends only on the choice of
the stable ruled surface $S = {\bf P}(F) \rightarrow C$.
Therefore a given stable ruled surface $S = {\bf P}(F) \in {\cal S}_C[3]$
defines uniquely its {\it associate Fano $3$-fold}  
$X =  X_S := M_C(2,K:3F)$, which is prime and of degree $16$.

\bigskip 

{\bf (5.6)} {\bf Lemma}. \   
{\sl
Let $F \rightarrow C$, $S = {\bf P}(F) \rightarrow C$,
$X = X_S$, $C_{\varepsilon}$, etc. be as in {\bf (5.4)}-{\bf (5.5)},
and let ${\Gamma}(X_S)$ be the curve of lines on $X_S$.
In addition, we shall assume that $X = X_S$ is general
(see {\bf (5.7)} below). 
Then the map
$$
\psi: Min(S) \rightarrow \ {\Gamma}(X_S),
\ C_{\varepsilon} \rightarrow l_{\varepsilon},
$$
\noindent
defined in {\bf (5.4)(b)}, is an isomorphism.
}

\bigskip

{\bf Proof}. \
Fix a general minimal section $C_{\varepsilon}$ of $S$,
and let {\bf (5.2)}(${\bf \ast}$) be the extension defined by
$C_{\varepsilon}$ as in {\bf (5.4)(b)}. 
By {\bf (5.2)} and {\bf (5.4)(c)}, the
sections $C_t \in Min \ {\bf P}(F) - \{ C_{\varepsilon} \}$
are in a $(1:1)$-correspondence with the elements of the family
\ $Sec_3^2({\Phi}_{\varepsilon}(C), [e_{\varepsilon}])$ \    
of $3$-secant planes ${\bf P}^2_t$ to the curve
${\Phi}_{\varepsilon}(C) \subset {\bf P}^4_{\varepsilon}$
which pass through the point $[e_{\varepsilon}]$. 
Clearly the projection from $[e_{\varepsilon}]$ sends
these planes ${\bf P}^2_t$ \ (1:1) \ to the $3$-secant
lines $L_t$ of $C^l_{\varepsilon} = C^l$, see {\bf (5.4)(d)}. 
In turn, the $3$-secant lines $L_t$ of $C^l_{\varepsilon}$
correspond, by {\bf (5.1.b)}, to the lines $l_t \subset X$
which do not intersect $l = l_{\varepsilon}$.
This way, the map 
$\psi: C_t \mapsto l_t$ is an isomorphism outside
the five lines $l_1,...,l_5$ which intersect
$l = l_{\varepsilon}$, see also {\bf (5.1.b)}.

The same local argument, but applied for
other general section $C_r$ of $S$
and the line $l_r \subset X$ corresponding
to $r$, implies that $C_t \rightarrow l_t$
is an isomorphism also outside 
the five $4$-secant lines $m_1,...,m_5$
of the line $m = l_r$ on $X$.
Since $l$ and $m$ are general, 
then the sets $\{ l_1,..,l_5 \}$
and $\{ m_1,..,m_5 \}$ are disjoint.
Therefore $\psi$ is an isomorphism.
{\bf q.e.d.}
 
\bigskip

{\bf (5.7)}
{\bf The associate ruled surface ${\bf S^X}$ of ${\bf X}$}. \
The inverse to {\bf (5.5)} is also true:
the general Fano $3$-fold $X = X_{16}$
is associate to a unique ruled surface $S$.
Indeed, let $C^l \subset {\bf P}^3_l$ be the curve defined
as in {\bf (5.1)} by the general line $l \subset X$,
and let $H = |K + \varepsilon|$ be the hyperplane system on
$C^l \subset {\bf P}^3_l$, where $K$ is the canonical class of $C^l$.
By {\bf (5.1)}{\bf (1)}-{\bf (2)}, for the general $l \subset X$, the curve
$C^l$ is a general curve of degree $7$ and of genus $3$
in ${\bf P}^3$. Therefore the complete linear system
$|H|$ sends $C^l$ isomorphically onto a curve
$C \subset {\bf P}^4$, and $C^l$ is a projection
of $C$ from a general point $[e] \in {\bf P}^4$.

By {\bf (5.2)}, the pair $(C,[e])$ represents
a minimal section $C_{\varepsilon}$ of the ruled surface
$S = {\bf P}(F) \in {\cal S}[3]$ defined by
an extension {\bf (5.2)}{\bf ($\ast$)} 
with $\varepsilon = \eta - \xi$. 
Clearly the surface $S$ is uniquely defined by
the pair $(X,l)$. By {\bf (5.4)}, 
the pair $(S,C_{\varepsilon})$ defines uniquely the pair
$(X,l)$, $l = l_{\varepsilon} \subset X$, 
and $X = X_S$ does not depend on the choice of the
minimal section $C_{\varepsilon}$ of $S$. 
This, and the isomorphism $\psi: Min(S) \rightarrow \Gamma(X_S)$,
$C_{\varepsilon} \mapsto l_{\varepsilon}$ from {\bf (5.6)}
(more precisely, the surjectivity of $\psi$),
imply the non-existence of other $S'$ with $X  = X_{S'}$.   

We call the unique ruled surface $S = S^X \rightarrow C$,
such that $X = X_S$, the {\it associate  ruled surface} of $X$;  
and call its base curve $C = C^X$ the {\it associate curve} of $X$. 
By \cite{CG}, the description {\bf (5.1)} of the double projection
imply that $C^X \cong C^l$ (see {\bf (5.4)(d)}) is the unique
curve such that the jacobian $J(C^X)$ is isomorphic, as a
principally polarized abelian variety, to the intermediate
jacobian $J(X)$ of $X$, see also \S 9 in \cite{M3}.

\bigskip
\bigskip 
\bigskip

\centerline{\large \bf \S 6. 
 The associate surface ${\bf S^X}$ is the Fano surface of ${\bf X}$}

\bigskip
\medskip

{\bf (6.1)}
{\bf The double projection from a line ${\bf l\subset\Sigma}$}. \
Let $l =  l_L \subset \Sigma$ be a line 
with an axis $L = {\bf P}(U_2^L)$ and a space
${\bf P}^3_L = {\bf P}(U_4^L)$, see {\bf (3.1)}.
By {\bf (2.4)(i)}, for any $u \in L$  
the tangent projective space ${\bf P}^6_u$
intersects $\Sigma$ along a cone $K_u \subset {\bf P}^6_u$
over the Veronese surface.
Let $W_L = \cup_{u \in L} \ {\bf P}^6_u$,
and let
$$
Z_L =  \cup_{u \in L} \ K_u = W_L \cap \Sigma .
$$
Then ${\bf P}^9_L := Span \ W_L = Span \ Z_L$
is a 9-dimensional subspace of ${\bf P}^{13} = Span \ \Sigma$. 
We denote by 
$\check{\bf P}^3_L = ({\bf P}^9_L)^{\perp} \subset \hat{\bf P}^{13}$ 
the space of linear equations of ${\bf P}^9_L \subset  {\bf P}^{13}$.

By its definition, the double projection
$$
{\pi}_{2.l}: \Sigma \rightarrow \check{\bf P}^3_L
$$
of $\Sigma \subset {\bf P}^{13}$ from $l$
is the restriction to $\Sigma$ of the projection
$\pi_{2.l}: {\bf P}^{13} \rightarrow \check{\bf P}^3_L$
from the subspace ${\bf P}^9_L = Span \ Z_L$. 
The codimension $2$ cycle $Z_L \subset \Sigma$
is the set of points $u \in \Sigma$
where ${\pi}_{2.l}$ is non-regular. Let also
$$
Z'_L = \{ u\in\Sigma : dim\ ({\bf P}^2_u \cap {\bf P}^3_L)\ge 1 \}.
$$
Since for $u \in \Sigma - Z'_L$ the lagrangian plane
${\bf P}^2_u$ intersects ${\bf P}^3_L$ at a point,
then the map
$$
\varphi_l: \Sigma - Z'_L \rightarrow {\bf P}^3_L, \
 u \mapsto {\bf P}^2_u \cap {\bf P}^3_L
$$
is regular.

\bigskip

{\bf (6.2)} {\bf Lemma}. \ 
{\sl
{\bf (i)} $Z'_L = Z_L$. 
\ {\bf (ii)}
There exists a natural identification
$\check{\bf P}^3_L \cong {\bf P}^3_L$
such that

\smallskip

\centerline{$\varphi_l = \pi_{2.l}|_{\Sigma-Z_L}$.} 
}

\bigskip

{\bf Proof of (i)}. \
Clearly $Z_L  \subset  Z'_L$;
and it rests to see that $Z'_L \subset Z_L$. 

Let $u \in Z'_L$, i.e. $dim \ ({\bf P}^2_u \cap {\bf P}^3_L) \ge 1$.
We have to see that $u \in Z_L$.

If ${\bf P}^2_u \subset {\bf P}^3_L$ then $u  \in l =  l_L$,
since the only lagrangian planes in ${\bf P}^3_L$ are the
planes ${\bf P}^2_w, w \in l = l_L$, see {\bf (3.1)}{\bf (ii)}.   

If $dim \ ({\bf P}^2_u \cap {\bf P}^3_L) =  1$
then the line $M = {\bf P}^2_u \cap {\bf P}^3_L$
will be an isotropic line in ${\bf P}^3_L$, 
since $M$ lies on the lagrangian plane ${\bf P}^2_u$.
Let $Q_L \subset G_L := G(1:{\bf P}^3_L)$
be the set of all isotropic lines in ${\bf P}^3_L$.
Clearly $Q_L$ contains the hyperplane section
${\sigma}_{10}(L) =
\{ M \subset {\bf P}^3_L: M \cap L \not= \emptyset \}$
= $\{$ the lines which lie in the lagrangian planes
${\bf P}^2_w, w \in L \}$.
Since the isotropic grassmannian $LG_2$ is a hyperplane
section of $G(2,6)$ (see {\bf (2.1)}), then
either $Q_L = G_L$, or $Q_L$ will coincide with the
hyperplane section ${\sigma}_{10}(L) \subset G_L$.
But if $Q_L = G_L = G(2,4)$,
then in ${\bf P}^3_L$ will lie two $3$-fold families of
lagrangian planes, which is impossible since the only lagrangian
planes in ${\bf P}^3_L$ are the planes ${\bf P}^2_w, w \in L$,
see {\bf (3.1)}{\bf (ii)}. 
Therefore $Q_L = {\sigma}_{10}(L)$. In particular 
$M \cap L \not= \emptyset$;
and since both $M$ and $L$ lie in ${\bf P}^3_L$
then the plane $Span(M \cup L)$ will be one,
say ${\bf P}^2_v$, of the planes
${\bf P}^2_w, w \in l  = l_L$, see {\bf (3.1)(i)}.
Since ${\bf P}^2_u$ intersects ${\bf P}^2_v$ along a line
then $u \in K_v$, see {\bf (2.4)(i)}. 
Therefore $Z'_L \subset \cup_{v \in L} \ K_v = Z_L$.
{\bf  q.e.d.}

\

{\bf Proof of (ii)}. \   
By the $Sp_3$-homogeneity of $LG_2$, it is enough to prove
{\bf (ii)} for a particular line $l \subset \Sigma$;   
and let $l = l_L  =  Span(e_{123},e_{423})$,
$L = {\bf P}(U_2^L) = {\bf P}(<e_2,e_3>)$,
$U_4^L = <e_1,e_2,e_3,e_4>$, etc.
be as in the proof of {\bf (3.1)}{\bf (ii)}. 
Then 
${\bf P}^3_L = {\bf P}(U_4^L) = {\bf P}^3(x_1:x_2:x_3:x_4)
\subset {\bf P}^5 = {\bf P}(V_6)$,
and
${\bf P}^9_L = Span \ Z_L = (x_{156} = x_{256} = x_{356} = x_{456} = 0)$
= $(y_{11} = y_{12} = y_{13} = z = 0) \subset {\bf P}^{13}$.
Therefore
$$
\check{\bf P}^3_L
= {\bf P}^3(x_{156}:x_{256}:x_{356}:x_{456})
\cong {\bf P}(U_4^L \otimes \wedge^2 \ V_6/U_4^L)
\cong {\bf P}^3_L .
$$

Under this identification, we shall verify first the coincidence
$\varphi_l =  \pi_{2l}|_{\Sigma-Z_L}$ over the open subset 
$$
\Sigma^o := \Sigma \cap (u=0) = exp(Sym^2 \ {\bf C}^3)
= \{ (1:X:\wedge^2 \ X : det \ X): {}^tX = X \}.
$$

By {\bf (2.1)}, the lagrangian 3-spaces
$U =  U_X \in \Sigma^o \subset G(3,6)^o$
are parameterized by the symmetric linear maps
$X = X_U: U_o \rightarrow  U_{\infty}$, such that
$U = <e_1 + X(e_1),e_2 + X(e_2),e_3+X(e_3)>$.
For such $U = U_X$ and $X = X_U$, the Pl\"ucker image
$pl(U_X) = {\bf P}(\wedge^3 \ U_X)$ is the same as the point
$exp(X) = (1: X : \wedge^2 \ X : \det \ X)$, ibid.
In other words, for $u \in \Sigma^o \subset {\bf P}^{13}$
one has: $u = exp(X_U)$,
where ${\bf P}(U) = {\bf P}^2_u$ is the lagrangian
projective plane of $u$.

\

Let $u \in \Sigma^o, u \not\in Z_L$, 
let ${\bf P}^2_u =  {\bf P}(U)$ be the lagrangian
plane of $u$, and let $X = X_U$.
We shall compute separately
$\varphi_l(u)$ and $\pi_{2.l}(u)$. 

On the one hand, $\varphi_l(u)$ is the intersection point
of ${\bf P}^2_u$ and ${\bf P}^3_L$. 
Since
${\bf P}^2_u = {\bf P}(U)
= {\bf P}(<e_1+X(e_1),e_2+X(e_2),e_3+X(e_3)>)$, 
and ${\bf P}^3_L = {\bf P}(U_4^L)
= {\bf P}(<e_1,e_2,e_3,e_4>)$,
then
$U \cap U_4^L =
<det(X_{11})e_1+det(X_{12})e_2+det(X_{13})e_3+det(X)e_4>$,
i.e.
$$
{\varphi}(u) = (det\ X_{11}:det\ X_{12}:det\ X_{13}:det\ X))
\in {\bf P}^3_L = {\bf P}^3(x_1:x_2:x_3:x_4),
$$
where for any $(i,j)$, $X_{ij}$ is the $(i,j)^{th}$ adjoint matrix of $X$. 

On the other hand,
as a point of $\Sigma \subset {\bf P}^{13}$,
$u = exp(X)  = (1:X:Y:z) \in \Sigma \subset {\bf P}^{13}$,
where $X = X_U$, $Y = \wedge^2 \ X$, $z = det \ X$. 
Since ${\pi}_{2.l}$ is the projection
from ${\bf P}^9_L$
onto $\check{\bf P}^3_L$
= ${\bf P}^3(x_{156}:x_{256}:x_{356}:x_{456})$
= ${\bf P}^3(y_{11}:y_{12}:y_{13}:z)$,  
then
$$
{\pi}_{2.l}(u) = {\pi}_{2l}(exp(X)) =
((\wedge^2\ X)_{11}:(\wedge^2\ X)_{12}:(\wedge^2\ X)_{13}:det\ X)
\in {\bf P}^3(y_{11}:y_{12}:y_{13}:z) .        
$$

Since $(\wedge^2\ X)_{ij} = det \ X_{ij}$, $i,j = 1,2,3$, \   
then ${\varphi}_l(u) = {\pi}_{2.l}(u)$, \ 
under the above identification
between ${\bf P}^3(x_1:x_2:x_3:x_4)$ = ${\bf P}^3_L$
and ${\bf P}^3(y_{11}:y_{12}:y_{13}:z)$ = $\check{\bf P}^3_L$.

The map ${\varphi}_l$ is evidently surjective onto
the $3$-space ${\bf P}^3_L$; and since 
${\varphi}_l$ coincides with the projective-linear map
${\pi}_{2l}$ over an open subset of $\Sigma$, 
then $\varphi_l$ is a restriction of a
projection of $\Sigma \subset {\bf P}^{13}$
from a $9$-space ${\bf P}^9_o \subset {\bf P}^{13}$.
Since ${\pi}_{2.l}$ is a projection from the $9$-space
${\bf P}^9_L = Span(Z_L)$, and since $Z_L$ is the non-regular
locus of ${\varphi}_l$, then
${\bf P}^9_o \supset Span(Z_L) = {\bf P}^9_L$.
Therefore ${\bf P}^9_o = {\bf P}^9_L$,
i.e. ${\varphi}_l = {\pi}_{2.l}|_{\Sigma - Z_L}$.
{\bf q.e.d.}

\bigskip 

{\bf (6.3)} {\bf The double projection from a line
${\bf l \subset X}$ and the ${\bf Sp_3}$-geometry}. \

Let $l \subset X = \Sigma \cap {\bf P}^{10}_X$ be a line,
let $L \subset {\bf P}^5$ be the axis of $l = l_L$,
and let ${\bf P}^3_L \subset {\bf P}^5$ be the space of $L$,
see {\bf (3.1)}.
The double projection ${\pi}_{2.l}: X \rightarrow {\bf P}^3_l$,
described in {\bf (5.1)}, 
is also the restriction to $X$ of the double projection
${\pi}_{2.l} :\Sigma \rightarrow {\bf P}^3 = {\bf P}^3_L$.
By {\bf (6.1)}, 
$\pi = {\pi}_{2.l}: X \rightarrow {\bf P}^3_L$ is non-regular
along the subset $z_l = Z_L \cap {\bf P}^{10}_X$,
and by {\bf (5.1)}, the set $z_l = l \cup l_1 \cup ...\cup l_e$,
where $l_1,...,l_e$ are the $e =  e(l) \le 5$ lines 
on $X$ which intersect $l$.   
By {\bf (6.2)}, for any $u \in X - z_l$ the double projection
${\pi} = {\pi}_{2.l}:X \rightarrow {\bf P}^3_l$,
described in {\bf (5.1)}, is given by:

\bigskip

${\bf (\ast)} \hspace{4.0cm}
\pi : u \mapsto {\pi}(u) = {\bf P}^2_u \cap {\bf P}^3_L$.    

\bigskip 

This identifies ${\bf P}^3_L$ and ${\bf P}^3_l$;  
and we shall find the curve $C^l \subset {\bf P}^3_L$,
see {\bf (5.1)(1)}. 

\

Let $l_1,...,l_e$, $e \le 5$ be the lines on $X$
which intersect $l$; and assume  for simplicity that
$e =  5$, see {\bf (5.1)(1)(iv)}. 
Let $L_1,...,L_5$ be the axes of $l_1,...,l_5$,
see {\bf (3.1)}.
Let $u_i  = l \cap l_i$, and let
${\bf P}^2_{u_i} \subset {\bf P}^5$ be the lagrangian
plane of $u_i$.
Since $u_i \in l$ then the axis $L$ of $l  = l_L$
lies in ${\bf P}^2_{u_i}$, and since $u_i \in l_i$
then $L_i \subset {\bf P}^2_{u_i}$.
Therefore $L_i$ intersects $L$, and let
$x_i = L \cap L_i$, $i = 1,...,5$.
The point $x_i = x(l+l_i)$ is the vertex of the
conic $q_i = l + l_i$, i.e. $x_i \in S_X$;
and let $\delta_L = \{ x_1,...,x_5 \}$.
Since $l_1,...,l_5$ are all the lines on $X$
which intersect $l$ then
$$
L \cap S_X = \delta_L.$$
In the general case $\delta_L$ is the intersection
$0$-cycle of $L$ and $S_L$, see {\bf (5.1)(1)(iv)}.  
Let
$$
C_L = \mbox{the closure of } S_X \cap ({\bf P}^3_L-L).  
$$

Then $S_L \cap {\bf P}^3_L =  C_L  \cup \delta_L$;
and we shall see that $C_L = C^l$.

\

Let first $x \in C_L \subset S_X  \cap {\bf P}^3_L$.
Since $x \in S_X$, then by {\bf (4.5.c)} there
exists a unique conic $q = q_x \subset X$ with vertex $x$;
and by $({\bf \ast})$ the double projection
${\pi}_{2.l}$ contracts the conic $q_x$ to the
point $x = x(q) \in {\bf P}^3_L = {\bf P}^3_l$;
in particular $x \in C^l$, see {\bf (5.1)(1)(ii)}.
Therefore $C_L \subset C^l$.

The same argument in the inverse direction
gives $C^l \subset C_L$. 

Notice that for the general $l$ the curve $C^l = C_L$
does not intersect the axis $L$ of the
line $l = l_L$, see {\bf (5.1.a)}.   
By {\bf (5.1)(1)(iv)}, the five axes $L_1,...,L_5$ 
are exactly the five $4$-secant lines to $C^l = C_L$, 
see also {\bf (5.1.a)} and {\bf (5.1.b)}. 

\bigskip

{\bf (6.4)} {\bf Corollary}. \
{\sl
In the notation of {\bf (6.3)}, \ $C_L = C^l$.
} 

\bigskip

{\bf (6.5)} \
{\bf Theorem}. \
{\sl
Let $X$ be a general prime Fano threefold of degree $16$.
Then the associate ruled surface $S^X$ of $X$
is isomorphic to the Fano surface ${\cal F}(X)$
of conics on $X$.
}

\bigskip

{\bf Proof}. \
Let $C^X$ and $C_X$ be correspondingly the associate
and the dual curve of $X$. 
By  {\bf (4.6)}, it is enough to prove that the ruled 
surface $p^X: S^X \rightarrow C^X$ is isomorphic
to the ruled vertex surface $p_X: S_X \rightarrow C_X$
of $X$.

Let $l = l_L \subset X$ be a line,
and let $C_{\varepsilon} \cong C^X$ be
the minimal section of $S^X$ such that
$l = l_{\varepsilon}$,
see {\bf (5.4)}, {\bf (5.6)}.
Let also $C^l \subset {\bf P}^3_l$
be the curve defined by the double
projection from $l$ as in {\bf (5.1)},
and let
$C^l_{\varepsilon} \subset {\bf P}^3_{\varepsilon}$
be the curve defined in {\bf (5.4)(c)}.

By {\bf (5.4)(c)-(5.4)(d)} \
$C^X \cong C^l$; and by {\bf (6.4)} \
$C^l \cong C_L$.
Since both $C^X$ and $C_X$ are
curves of genus $3$, then $C_L \subset S_X$
is a section of $S_X$.  Therefore 

\medskip

\centerline{$C^X \cong C_L \cong C_X =: C .$}

\medskip

Let   
$Min'(S_X) =  \{ C_L : l = l_L \in {\Gamma}(X)\}$
be the family of all sections $C_L \subset S_X$. 
One can see separately that $Min'(S_X) = Min(S_X)$,
but in the proof we will use only the isomorphism  
$Min'(S_X) \cong Min(S^X)$. 
%
%
In fact, by {\bf (5.4)(c)-(d)} and {\bf (6.4)}
$$
C^l  = C^l_{\varepsilon} = C_L  \ \mbox{  in  } \ 
{\bf P}^3_l = {\bf P}^3_{\varepsilon} =  {\bf P}^3_L ;
$$
and then by {\bf (5.6)} the maps   

\begin{picture}(130,15)
\put(25,6){\makebox(0,0){$Min'(S_X)$}}
\put(55,6){\makebox(0,0){$Min(S^X)$}}
\put(90,6){\makebox(0,0){$\Gamma(X) =: \Gamma$ ,}}
\put(120,6){\makebox(0,0){$C_L \mapsto
         C_{\varepsilon} \mapsto l_{\varepsilon}$}}
\put(35,6){\vector(1,0){10.5}}
\put(65,6){\vector(1,0){12.5}}
\put(40,9){\makebox(0,0){$\xi$}}
\put(70,9){\makebox(0,0){$\psi$}}
\end{picture}

\noindent are isomorphisms. 

It rests to see that $S^X \cong S_X$.
For this, define the following
maps:

\medskip

\centerline{$\alpha: C \times \Gamma \rightarrow S_X, \
             \alpha : (c,l) \mapsto x.(c) ,$}

\medskip

\noindent 
where $l = l_L$, and
$x.(c) := (p_X)^{-1}(c) \cap C_L \in S_X$,
for $c \in C = C_X$; and  

\medskip

\centerline{$\beta : C \times \Gamma \rightarrow S^X, \
              \beta : (c,l)  \mapsto x^.(c) ,$} 

\medskip

\noindent
where $l = l_{\varepsilon}$, and
$x^.(c) := (p^X)^{-1}(c) \cap C_{\varepsilon} \in S^X$,
for $c \in C = C^X$.

In order to prove that $S_X \cong S^X$, it is enough 
to identify the fibers $\alpha^{-1}(\alpha(c,l))$
and ${\beta}^{-1}(\beta(c,l))$,
for any $(c,l) \in C \times \Gamma$. 

Let $(c,l) \in C \times \Gamma$. By the definition of $\alpha$ 

\medskip

\centerline{$\alpha^{-1}(\alpha(c,l)) \cong
       \{ C_{L_i} \subset S_X:
               l_{L_i} \subset X \ \& \ x.(c) \in C_{L_i} \},$}

\medskip

\noindent
and one of these sections is $C_L$, $l = l_L$.
Let $q_{x.(c)}$ be the unique conic on $X$
with vertex $x.(c)$, see {\bf (4.5.c)} or the proof of {\bf (4.4)}.
Since

\medskip

\centerline{$q_{x.(c)} = Q_{x.(c)} \cap {\bf  P}^{10}_X 
            = \{ u \in X : x.(c) \in {\bf P}^2_u \}$}

\medskip

\noindent
(see e.g. {\bf (6.3)} or the proof of {\bf (4.4)}),
then the line $l_i = l_{L_i} \subset X$ passes through $x.(c)$
\ {\it iff } \ $l_i \cap q_{x.(c)} \not= \emptyset$.
Therefore

\medskip

\centerline{$\alpha^{-1}(\alpha(c,l)) \cong
\{ l_i \in {\Gamma}(X): q_{x.(c)} \cap l_i \not= \emptyset \}.$}

\medskip

\noindent
Clearly, the line $l$ is one of these lines.

By the definition of $\beta$

\medskip

\centerline{$\beta^{-1}(\beta(c,l)) \cong
             \{ C_t \subset S^X : x^.(c) \in C_t \},$}

\medskip

\noindent 
and one of these sections is
$C_{\varepsilon}$, $l = l_{\varepsilon}$.
By {\bf (5.2)} and {\bf (5.4)(c)} 
(see also the proof of {\bf (5.6)})
the minimal sections $C_t$ of $S^X$, different from
$C_{\varepsilon}$, are parameterized
by the $3$-secant lines $M$ to
$C^l_{\varepsilon} \subset {\bf P}^3_{\varepsilon}$.
Since $C^l_{\varepsilon}
= C_L \subset {\bf P}^3_L =  {\bf P}^3_{\varepsilon}$,
then $C_t$ intersects $C_{\varepsilon}$ at the point $x^.(c)$
\ {\it iff} \ $x.(c) \in M$.
Therefore

\medskip

\centerline{$\beta^{-1}(\beta(c,l)) \cong
  \{ M \in Sec_3(C_L): x.(c) \in M \} \cup \{ C_{\varepsilon} \}.$} 

\medskip

By the proof of {\bf (5.1.b)} the $3$-secant lines  
$M$ to $C_L  = C^l$ are the proper ${\pi}$-transforms
of the lines $m \subset X$. 
Moreover, by {\bf (5.1)}{\bf (1)}{\bf (ii)} the points $x = x_q \in C^l$
are the same as the blowed-down
strict transforms $q^+ \subset X^+$ of the conics
$q \subset X$ which intersect $l$.
Therefore there exists a unique conic
$q = q_o \subset X$ such that $x_{q_o} = x.(c)$. 
Since by {\bf (6.2)-(6.3)(${\bf \ast}$)} the double projection 
${\pi} = {\pi}_{2.l}$ contracts $q_o$ to its vertex
$x_{q_o} =  x.(c)$ then $q_o = q_{x.(c)}$, and 

\medskip

\centerline{$\beta^{-1}(\beta(c,l)) \cong
    \{ l_i \in {\Gamma}(X): q_{x.(c)} \cap l_i \not= \emptyset \}
    \cong \alpha^{-1}(\alpha(c,l)).$} 

\medskip

\noindent
{\bf q.e.d.}

\bigskip

{\bf (6.6)} 
{\bf Remark}. \   
In fact, the coincidence between the fibers of
$\alpha$ and $\beta$ implies the isomorphsm between
the images ${\alpha}(C \times \Gamma) \subset S_X$
and ${\beta}(C \times \Gamma) \subset S^X$.
But $\alpha$ and $\beta$ are certainly surjective
since a conic $q \subset X$ always intersects 
the surface $R_X \subset X$ swept out by the
lines $l \subset X$. The general
conic $q$ on the general $X$ intersects exactly
$8$ lines on $X$ since $R_X \in |{\cal O}_X(4)|$,
see (6.4) case (v) in \cite{I1}. 
The virtual number $e(q) = 8$ of lines intersecting
a conic $q \subset X_{16}$ can be computed 
also by the Mori theory -- see \S (2.5)-(2.6) and
table (2.8.2) in \cite{T}.  

On the general ruled surface $S = S^X \rightarrow C^X  = C$
of invariant $e(S) = 3$, 
one can compute the same number $e(x) = 8$ of minimal
sections through its general point $x$ as follows
(see {\bf (5.2)}): \ If $x \in S$ is general then
the elementary transformation $elm_x: S \rightarrow S_x$
sends $S$ to a general $S_x \in {\cal S}_C[2]$.
The proper images of the minimal sections of $S$ through $x$
are all the minimal sections of $S_x$.
$S_x$ always has minimal sections, and let $C_{o} \subset S_x$
be one of them. Since $S_x$ is general, then 
$\varepsilon_o = C_o^2 := {\cal O}_{C_o}(C_o) \in Pic^2(C_o)$
= $Pic^2(C)$ is general, and the linear system
$|K_C + \varepsilon_o|$ sends $C$ isomorphically to a
space curve $C^3_6 \subset {\bf P}^3$ of genus $3$ and
of degree $6$. The minimal sections of $S_x$,
different from $C_o$, are in a 1:1 correspondence with
the bisecant lines to $C^3_6$ from the (general)
extension-class point $[e_{\varepsilon_o}] \in {\bf P}^3$
defined by the section $C_o$. These lines are $7$,
since the projection of $C^3_6 \subset {\bf P}^3$
from a general point in ${\bf P}^3$ is a plane sextic
with $7$ double points.


\bigskip
 

%
{\small
}

\bigskip

{\small ${}$\hspace{1.2cm}{Atanas Iliev}

${}$\hspace{1.2cm}Institute of Mathematics

${}$\hspace{1.2cm}Bulgarian Academy of Sciences

${}$\hspace{1.2cm}Acad. G. Bonchev Str., 8

${}$\hspace{1.2cm}1113 Sofia, Bulgaria

${}$\hspace{1.2cm}e-mail: ailiev@math.bas.bg
\hspace{2.7cm}
}

\end{document}